\documentclass[reqno]{amsart}
\usepackage{amsmath}
\usepackage{amsfonts}
\usepackage{amssymb,enumerate}
\usepackage{amsthm}
\usepackage[all]{xy}
\usepackage{rotating}
\usepackage{hyperref}
\usepackage{color}
\theoremstyle{plain}
\newtheorem{lem}{Lemma}[section]

\newtheorem{prop}[lem]{Proposition}
\newtheorem{thm}[lem]{Theorem}
\newtheorem{conj}[lem]{Conjecture}

\theoremstyle{definition}

\newtheorem{question}[lem]{Question}

\newtheorem{para}[lem]{}

\newcommand{\cat}[1]{\mathcal{#1}}

\newcommand{\catd}{\cat{D}}

\newcommand{\ges}{\geqslant}

\newcommand{\pd}{\operatorname{pd}}

\newcommand{\id}{\operatorname{id}}	
\newcommand{\fd}{\operatorname{fd}}

\newcommand{\depth}{\operatorname{depth}}	
\newcommand{\rank}{\operatorname{rank}}	

\newcommand{\edim}{\operatorname{edim}}

\newcommand{\len}{\operatorname{len}}
\newcommand{\type}{\operatorname{type}}

\newcommand{\HH}{\operatorname{H}}
\newcommand{\Hom}{\operatorname{Hom}}

\newcommand{\s}{\mathfrak{S}}

\newcommand{\im}{\operatorname{Im}}
\newcommand{\shift}{\mathsf{\Sigma}}

\newcommand{\Ker}{\operatorname{Ker}}

\newcommand{\ideal}[1]{\mathfrak{#1}}

\newcommand{\fm}{\ideal{m}}
\newcommand{\fn}{\ideal{n}}
\newcommand{\fp}{\ideal{p}}
\newcommand{\fq}{\ideal{q}}

\newcommand{\xra}{\xrightarrow}
\newcommand{\xla}{\xleftarrow}

\newcommand{\ecodepth}{\operatorname{ecodepth}}

\renewcommand{\geq}{\geqslant}
\renewcommand{\leq}{\leqslant}
\renewcommand{\ker}{\Ker}

\newcommand{\Ext}[4][R]{\operatorname{Ext}_{#1}^{#2}(#3,#4)}	
\newcommand{\Rhom}[3][R]{\mathbf{R}\!\operatorname{Hom}_{#1}(#2,#3)}	

\newcommand{\Otimes}[3][R]{#2\otimes_{#1}#3}
\renewcommand{\Hom}[3][R]{\operatorname{Hom}_{#1}(#2,#3)}	
\newcommand{\Tor}[4][R]{\operatorname{Tor}^{#1}_{#2}(#3,#4)}

\newcommand\mult{\operatorname{mult}}

\newcommand{\od}{\operatorname{Mod}}

\newcommand{\ggl}{\operatorname{GL}}
\newcommand{\ta}{\operatorname{\textsf{T}}}

\newcommand{\yext}[4][R]{\operatorname{YExt}_{#1}^{#2}(#3,#4)}

\newcommand{\Taglm}{\ta^{\ggl(W)_0\cdot M}}

\newcommand{\pf}{\operatorname{pf}}
\newcommand{\DD}{\operatorname{D}}
\newcommand{\LD}{\operatorname{L}}

\numberwithin{equation}{lem}

\begin{document}

\bibliographystyle{amsplain}

\author{Saeed Nasseh}
\address{Department of Mathematical Sciences\\
Georgia Southern University\\
Statesboro, GA 30460, USA}
\email{snasseh@georgiasouthern.edu}
\urladdr{https://sites.google.com/site/saeednasseh/home}

\author{Sean K. Sather-Wagstaff}
\address{School of Mathematical and Statistical Sciences,
Clemson University,
O-110 Martin Hall, Box 340975, Clemson, S.C. 29634
USA}
\email{ssather@clemson.edu}
\urladdr{https://ssather.people.clemson.edu/}

%\thanks{
%This material is based on work supported by North Dakota EPSCoR and
%National Science Foundation Grant EPS-0814442.
%Sean Sather-Wagstaff was supported in part by a grant from the NSA}

\title[Applications of DG Algebra Techniques in Commutative Algebra]{Applications of Differential Graded Algebra Techniques in Commutative Algebra}

%\date{\today}

\dedicatory{To David Eisenbud on the occasion of his 75th birthday}

%\keywords{Differential graded algebra, DG algebra, }

\begin{abstract}
Differential graded (DG) algebras are powerful tools from rational homotopy theory. We survey some recent applications of these in the realm of homological commutative algebra.
\end{abstract}

\maketitle

\tableofcontents

\section{Introduction} \label{sec0}
Throughout this paper, the term ``ring'' is short for ``commutative noetherian ring with identity.''

In algebraic topology, it is incredibly useful to know that the singular cohomology of a manifold has a natural algebra structure. Similarly, in commutative algebra the fact that certain Ext and Tor modules carry algebra structures is a powerful tool. Both of these notions arise by considering differential graded (DG) algebra structures on certain chain complexes. In short, a DG algebra is a chain complex that is also a graded commutative ring, where the differential and multiplication are compatible; see Section~\ref{sec2} for definitions and background material.

Avramov, Buchsbaum, Eisenbud, Foxby, Halperin, Kustin, and others pioneered the use of DG algebra techniques in homological commutative algebra. The idea is to prove results about rings by broadening one's context to include vast generalizations. A deep, rich sample of the theory and applications can be found in Avramov's lecture notes~\cite{avramov:ifr}.

The current paper is a modest follow-up to {\it op.\ cit.}, documenting a few applications that have appeared in the twenty-some years since {\it op.\ cit.} appeared. To be clear, we focus on applications: results whose statements make no reference to DG algebras but whose proofs use them extensively. Furthermore, this survey is by no means comprehensive. We focus on small number of some of our favorite applications, limited by constraints of time and space.

Most of the sections below begin by describing an application with little reference to DG algebras. This is followed by a certain amount of DG background material, but generally only enough to give a taste for the material. The sections conclude with an indication of how the DG technology helps to obtain the application.

As we noted above, the work of David Eisenbud is foundational in this area, especially the paper~\cite{buchsbaum:asffr} with Buchsbaum; see~\ref{para20201006c}. Those of us working in this area owe him a huge debt of gratitude for this and other seminal work in the field.

\section{Growth of Bass and Betti numbers}\label{sec2}

In this section, let $(R,\fm,k)$ be a local ring with $d:=\depth R$. The \emph{embedding codepth} of $R$, denoted $c:=\ecodepth R$, is defined to be $e-d$, where $e:=\edim R$ is the minimal number of generators of $\fm$. Cohen's Structure Theorem states that the $\fm$-adic completion $\widehat{R}$ admits a minimal Cohen presentation, i.e., there is a complete
regular local ring $(P,\fp, k)$ and an ideal $I\subseteq \frak p^2$
such that $\widehat{R}\cong P/I$. Note that the projective dimension $\pd_P(\widehat{R})$, i.e., the length of the minimal free resolution of $\widehat{R}$ over $P$, is equal to $\ecodepth R$.

Fundamental invariants of a finitely generated $R$-module $M$ are the Bass and Betti numbers. These numerically encode structural information about the module $M$, e.g., the minimal number of generators and relations and higher degree versions of these. A hot topic of research in commutative algebra is the growth of these sequences. In this section, we describe some recent progress by Avramov on this subject including how he uses DG techniques to get information about these invariants. Along the way, we also present foundational material about the DG~context.

\subsection*{Bass numbers, Betti numbers, and a question of Huneke}\

%\begin{para}\label{para20200901s}
Let $M$ be a finitely generated $R$-module. For each integer $i$, the $i$th \emph{Bass number} and the $i$th \emph{Betti number} of $M$ are defined to be, respectively
\begin{gather}
\mu_R^i(M):=\rank_k\left(\Ext[R]{i}kM\right)\notag\\ \beta^R_i(M):=\rank_k\left(\Ext[R]{i}Mk\right)=\rank_k\left(\Tor[R]{i}Mk\right).\notag
\end{gather}
The \emph{Bass series} and the \emph{Poincar\'{e} series} of $M$ are the formal power series
\begin{equation}\label{eq20200902a}
I^M_R(t):=\sum_{i\in \mathbb{Z}}\mu_R^i(M) t^i\quad\quad\quad P_M^R(t):=\sum_{i\in \mathbb{Z}}\beta^R_i(M) t^i.
\end{equation}
In case that $M=R$, the Bass numbers and Bass series are denoted $\mu_R^i$ and $I_R(t)$.
%\end{para}

In this section we are concerned with the following unpublished question of Huneke; see~\cite{sather:bnsc}. This question is motivated in part by the fact that $R$ is Gorenstein if and only if its Bass numbers are eventually $0$.

\begin{question}\label{question20200813a}
Let $R$ be a Cohen-Macaulay local ring. If $\{\mu_R^i\}$ is bounded, must $R$ be Gorenstein? If $\{\mu_R^i\}$ is bounded above by a polynomial in $i$, must $R$ be Gorenstein? If $R$ is not Gorenstein, must $\{\mu_R^i\}$ grow exponentially?
\end{question}

Very little progress has been made on this question. Christensen, Striuli, and Veliche~\cite{christensen:gmirlr} conduct a careful analysis of several special cases of this and other related questions. Other progress comes from Jorgensen and Leuschke~\cite{jorgensen:gbscm} and Borna, Sather-Wagstaff, and Yassemi~\cite{borna:vhgbscm, sather:bnsc}.

In this section, we focus on work of Avramov~\cite{avramov:cslrec3} on this question for non-Gorenstein rings $R$ with $c=\ecodepth(R)\leq 3$ which includes the following result. The proof relies heavily on DG techniques as we explain in the next subsections. This is a true application of DG tools, as the statement makes no mention of DG algebras, though they are used extensively in the proof.

\begin{thm}[\protect{\cite[Theorem 4.1]{avramov:cslrec3}}]\label{thm20200718b}
If $c\leq 3$ and $R$ is not Gorenstein, then there is a real number $\gamma_R>1$ such that for all $i\geq 1$
\begin{equation}\label{eq20200718c}
\mu_R^{d+i}\geq \gamma_R \mu_R^{d+i-1}
\end{equation}
with two exceptions for $i=2$: If $I=(wx,wy)$ or $I=(wx,wy,z)$, where $x,y\in P$ is a regular sequence, $w\in P$, and $z\in \fp^2$ is a $P/(wx,wy)$-regular element, then $\mu_R^{d+2}=\mu_R^{d+1}=2$. If $R$ is Cohen-Macaulay, the inequality~\eqref{eq20200718c} holds for all $i$.
\end{thm}

\subsection*{DG algebra resolutions and DG modules}\

Let $S$ be a ring. A \emph{associative, commutative differential graded $S$-algebra}
(\emph{DG $S$-algebra} for short) is a chain complex
$A=\cdots\to A_2\xra{\partial^A_2}A_1\xra{\partial^A_1}A_0\to 0$ such that $A^{\natural}=\bigoplus_{i\geq 0}A_i$ has the structure of a graded commutative $S$-algebra
\begin{enumerate}[$\bullet$]
\item
for all $a,b\in A$ the equality $ab=(-1)^{|a||b|}ba$ holds, and $a^2=0$ if the homological degree $|a|$ is odd
\end{enumerate}
that satisfies the \emph{Leibniz rule}
\begin{enumerate}[$\bullet$]
\item
for all $a,b\in A$ we have $\partial^A(ab)=\partial^A(a)b+(-1)^{|a|}a\partial^Ab$, i.e., the assignment $a\otimes b\mapsto ab$ describes a chain map $\Otimes[S]AA\to A$.
\end{enumerate}
The DG $S$-algebra $A$ is called \emph{homologically degreewise noetherian} if $\HH_0(A)$ is noetherian, and each $\HH_0(A)$-module $\HH_i(A)$ is finitely generated.

Examples of homologically degreewise noetherian DG $S$-algebras include $S$ itself, considered as a complex concentrated in degree $0$, and the Koszul complex $K^S(\mathbf{\underline{x}})$ over $S$ on a sequence $\mathbf{\underline{x}}=x_1,\ldots,x_n$ in $S$ with the exterior algebra structure.

A \emph{morphism} of DG $S$-algebras is a chain map $f\colon A\to B$ such that for all $a,a'\in A$ we have $f(aa')=f(a)f(a')$ and $f(1)=1$. A \emph{quasiisomorphism} of DG algebras is a morphism that is a quasiisomorphism, i.e., such that the induced map on homology is an isomorphism in each degree.
A \emph{DG algebra resolution of an $S$-algebra $T$} is a quasiisomorphism $F\xra{\simeq} T$ of DG $S$-algebras such that
each $F_i$ is free over $S$. Several examples of DG algebra resolutions are given below starting with~\ref{para20200803y}.

In case that $(S,\fn)$ is a local ring, a DG $S$-algebra $A$ is called \emph{local} if it is homologically degreewise noetherian and $\HH_0(A)$ is a local $S$-algebra. In this case, setting $\fn_0$ to be the preimage of $\fm_{\HH_0(A)}$ in $A_0$, we let $\fm_A=\cdots\xra{\partial^A_2}A_1\xra{\partial^A_1} \fn_0\to 0$ be the \emph{augmentation ideal} of $A$. By definition it is a subcomplex of $A$. Moreover, it is a \emph{DG ideal} of $A$ meaning that it absorbs multiplication by elements of $A$. In this situation, we say that $(A,\fm_A,A/\fm_A)$ is a local DG $S$-algebra. As an example, if $\mathbf{\underline{x}}\in \fn$ is a sequence of $n$ elements, then $K=K^S(\mathbf{\underline{x}})$ is a local DG $S$-algebra with the augmentation ideal $\fm_{K}=0\to S\to\cdots\to S^n\to \fn\to 0$ and $K/\fm_K\cong S/\fn$.

\begin{para}\label{para20201006a}
A construction of Tate~\cite{tate:hnrlr} (see also Avramov~\cite[Proposition 2.2.8]{avramov:ifr}) guarantees the existence of a DG algebra resolution $F$ of $\widehat{R}$ over $P$, where each $F_i$ is finitely generated and free over $P$ and
$F_i=0$ for all $i>\pd_P(\widehat{R})$.
\end{para}

Examples of DG algebra resolutions include the following.

\begin{para}\label{para20200803y}
If $I$ is generated by a $P$-regular sequence (that is, if $R$ is a formal complete intersection), then the Koszul complex $K^P(I)$ on a minimal generating sequence for $I$ is a DG algebra resolution of $\widehat{R}$ over $P$.
\end{para}

\begin{para}\label{para20201006b}
If $\pd_P(\widehat{R})=2$, then it follows from the Hilbert-Burch Theorem~\cite[Theorem~20.15]{eisenbud:ca} that there is an element $f\in P$ and a matrix $Y$ of size $n\times (n-1)$ such that the minimal $P$-free resolution of $\widehat{R}$ can be chosen with the form
\begin{equation*}%\label{eq20200813d}
0\to P^{\oplus n-1}\xra{Y} P^{\oplus n}\xra{X} P\to \widehat{R}\to 0
\end{equation*}
with $X=f\left(\det(Y_1),\ldots,(-1)^{j-1}\det(Y_j),\ldots,(-1)^{n-1}\det(Y_n)\right)$, where $Y_j$ is the minor obtained from $Y$ by deleting the $j$-th row. Herzog~\cite{herzog:kadla} describes a DG algebra structure on this resolution, as follows. Let $\{a_1,\ldots,a_n\}$ be a basis for $P^{\oplus n}$ and $\{b_1,\ldots,b_{n-1}\}$ be a basis for $P^{\oplus n-1}$, and set
\begin{gather}
(a_i)^2=0\notag\\
a_i\cdot a_j=-a_j\cdot a_i=\sum_{t=1}^{n-1}(-1)^{i+j+t+1}\det(Y_{ij,t})fb_t\qquad \text{for}\ i<j\notag
\end{gather}
where $Y_{ij,t}$ denotes the minor obtained from $Y$ by deleting rows $i,j$ and column $t$.
\end{para}

\begin{para}\label{para20201006c}
Assume $\pd_P(\widehat{R})=3$. Buchsbaum and Eisenbud~\cite{buchsbaum:asffr} show that the minimal free resolution of $\widehat{R}$ over $P$ has the structure of a DG algebra, though the explicit structure of the resolution is not given.

Let $n$ denote the minimal number of generators for $I$, and assume that $R$ is  Gorenstein (that is, $I$ is a Gorenstein ideal). Then it is shown in {\it op.\ cit.} that the minimal free resolution of $\widehat{R}$ over $P$ is of the form
\begin{equation}\label{eq20200814a}
0\to P\xra{Z} P^{\oplus n}\xra{Y} P^{\oplus n}\xra{X} P\to \widehat{R}\to 0
\end{equation}
for some $n\times n$ alternating matrix $Y$ with entries in $\fp$ and $$X=\left(\pf(Y_{1}),\ldots,(-1)^{j-1}\pf(Y_{j}),\ldots,(-1)^{n-1}\pf(Y_{n})\right)$$ where $Y_{i}$ denotes the matrix obtained from $Y$ by deleting row $i$ and column $i$.  (Here, $\pf$ is the Pfaffian; see~\cite{bruns:cmr} for details.) Also, $Z=\Hom[P]XP$.

An explicit DG algebra structure on~\eqref{eq20200814a} is given by Avramov~\cite{MR485906} as follows. Let $\{a_1,\ldots,a_n\}$ be a basis for $P^{\oplus n}$ in degree $1$, let $\{b_1,\ldots,b_{n}\}$ be a basis for $P^{\oplus n}$ in degree $2$, and let $\{c\}$ be a basis for $P$ in degree $3$. Define
\begin{gather}
(a_i)^2=0\qquad\qquad
a_i\cdot b_j=b_j\cdot a_i=\delta_{ij}c\notag\\
a_i\cdot a_j=-a_j\cdot a_i=\sum_{t=1}^{n}(-1)^{i+j+t}\rho_{ijt}\pf(Y_{ijt})b_t\qquad \text{for}\ i<j\notag
\end{gather}
where $Y_{ijt}$ is the matrix obtained from $Y$ by deleting rows $i,j,t$ and columns $i,j,t$, and $\delta_{ij}$ is the Kronecker delta, and $$\rho_{ijt}=\begin{cases}-1& i<t<j\\ 1 &\text{otherwise.}\end{cases}$$
\end{para}

There are many examples of DG algebra resolutions in the monomial situation. We summarize a few here very briefly and point the reader to references for more details. Let $S=k[x_1,\ldots,x_n]$ be a polynomial ring over a field $k$, and let $I$ be a monomial ideal in $S$, i.e., an ideal generated by monomials. A general DG algebra resolution of $S/I$ is given by Taylor~\cite{taylor:igmrs} and Bayer, Peeva, and Sturmfels~\cite{MR1618363}, but it is not minimal in general. In the following cases the minimal free resolution of $S/I$ over $S$ has a DG algebra structure: stable ideals (see Eliahou and Kervaire~\cite{MR1037391} or Peeva~\cite{MR1407879}), matroidal ideals (see Sk\"{o}ldberg~\cite{skoldberg1}), and ideals of the form $I=fJ$, where $J$ is a monomial ideal in $S$ and $f$ is the least common multiple of the generators of $J$ (see Katth\"{a}n~\cite{MR3862675}).

\begin{para}\label{para20200718m}
It is important to note for~\ref{para20201006a} that the minimal free resolution of $\widehat{R}$ over $P$ may not admit a DG algebra structure in general. Examples for this are given by Khinich (as documented in~\cite{MR0349816}), Avramov~\cite{MR601460}, and Katth\"{a}n~\cite{MR3862675}. The example of Katth\"{a}n is generic and disproves a claim by Bayer, Peeva, and Sturmfels~\cite{MR1618363}.
\end{para}

\begin{para}\label{para20200718a}
In contrast to~\ref{para20200718m}, if $R$ satisfies one of the following conditions, then the minimal free resolution of $\widehat{R}$ over $P$  admits a DG algebra structure:
\begin{enumerate}[\rm(a)]
\item
$c\leq 3$, by~\ref{para20201006b}--\ref{para20201006c};
\item
$c=4$ and $R$ is Gorenstein, by Kustin and Miller~\cite{kustin:gacfct, kustin:asmrgrecf};
\item
$c=4$, $R$ is Cohen-Macaulay, almost complete intersection, and $1/2\in R$, by Kustin~\cite{MR1292771};
\item
$R$ is complete intersection, by~\ref{para20200803y};
\item
$R$ is one link from a complete intersection, by Avramov, Kustin, and Miller~\cite{avramov:psmlrsec};
\item
$R$ is two links from a complete intersection and is Gorenstein, by Kustin and Miller~\cite{MR713381}.
\end{enumerate}
\end{para}

The translation to DG algebras uses the following DG analogue of modules.

\begin{para}\label{para20201021a}
Let $B$ be a DG $S$-algebra. A \emph{DG $B$-module} is an $S$-complex $M$ such that $M^{\natural}=\bigoplus_{i}M_i$ is a graded $A^{\natural}$-module satisfying the Leibniz rule.
For the case of $S$ considered as a DG $S$-algebra, the DG $S$-modules are just the $S$-complexes. A DG $B$-module $M$ is \emph{homologically bounded} if $\HH_i(M)=0$ for all $|i|\gg 0$; it is \emph{homologically finite} if $\bigoplus_i\HH_i(M)$ is a finitely generated $\HH_0(B)$-module.

Let $M$ be a DG $B$-module. The \emph{trivial extension} $B\ltimes M$ is the DG algebra with the underlying complex $B\oplus M$ equipped with the product that is given as follows:
$$
(b,m)(b',m'):=(bb',bm'+(-1)^{|m||b'|}b'm).
$$
\end{para}

\subsection*{Avramov's machine}\

Some of our favorite applications of DG techniques use the following tool which Kustin~\cite{MR1132435} calls \emph{Avramov's machine}.

\begin{para}\label{para20200803g}
Let $\mathbf{x}$ be a minimal generating sequence
for $\fm$, and let $\mathbf{y}$
be a minimal generating sequence for the maximal ideal $\fp$.
Since $P$ is a regular local ring, the Koszul complex $K^P(\mathbf{y})$ is a minimal free resolution
of $k$ over $P$.
Since $K^P(\mathbf{y})\simeq k$, we obtain the following diagram of DG algebra quasiisomorphisms:
\begin{equation}\label{eq201208c}
K^{R}(\mathbf{x})\xra{\simeq}K^{\widehat{R}}(\mathbf{x}\widehat{R})\xla{\cong}\Otimes[P]{K^P(\mathbf{y})}{\widehat{R}}\xla\simeq\Otimes[P]{K^P(\mathbf{y})}{F}\xra\simeq \Otimes[P]{k}{F}=:A.
\end{equation}
The assumptions on $F$ in~\ref{para20201006a} imply that $A$ is a finite-dimensional DG
$k$-algebra.
It follows that $\Tor[P]{}{\widehat{R}}{k}$ inherits the structure of a finite-dimensional DG $k$-algebra; this is the \emph{Tor algebra}. If $F$ is minimal, e.g., in any of the cases from~\ref{para20200718a}, the algebra $A$ has zero differential, so
$A\cong\Tor[P]{}{\widehat{R}}{k}$.
\end{para}

Rationality of Poincar\'{e} series, which we discuss next, is an important application of Avramov's machine.

\begin{para}\label{para20200809s}
Consider the notation from~\ref{para20200803g}. In this paragraph, assume that one of the conditions (a), (b), (e), or (f) in~\ref{para20200718a} holds. Using the fact that the minimal free resolution of $\widehat{R}$ over $P$ has a DG algebra structure, Avramov, Kustin, and Miller~\cite{avramov:psmlrsec} give a factorization $P\xra{\varphi}Q\xra{\psi}\widehat{R}$
of the canonical map $P\to \widehat{R}$ such that $\varphi$ is complete intersection and $\psi$ is Golod (see, e.g., \cite{MR846439} for the definition). Then they invoke a result of Levin~\cite{Levin} to conclude the following;

\begin{enumerate}[($*$)]
\item
the Poincar\'{e} series of every finitely generated $R$-module is rational with common denominator.
\end{enumerate}
In case (d) of~\ref{para20200718a} where $R$ is a complete intersection, conclusion ($*$) was proved for $P_k^R(t)$ by Tate~\cite{tate:hnrlr} and, in general, by Gulliksen~\cite{MR364232} and Avramov~\cite{avramov:vpd}. In case~(c) of~\ref{para20200718a}, conclusion ($*$) was proved by Kustin and Palmer~\cite{MR1295961}.
\end{para}

\subsection*{Growth rates in embedding codepth at most $3$}\

%As we mentioned earlier, in \cite{avramov:cslrec3}, Avramov specifically focuses on the case where $c\leq 3$ and investigates the structure of such rings along with the Poincar\'{e} and Bass series $P^R_k(t)$ and $I_R(t)$. In this subsection, we indicate how Avramov's machine contributes to the proof.

With these tools in hand, the proof of Theorem~\ref{thm20200718b} proceeds in the following steps. First, consider the following structure result for the Tor algebra.

\begin{para}\label{para20200816z}
Assume that $c\leq 3$. Using the notation from~\ref{para20200803g}, we know that $A$  is a finite-dimensional DG algebra with zero differential. In this case, by~\cite{avramov:psmlrsec} and~\cite{weyman:sfrl3}, the ring $R$ belongs to one of the following classes\vspace{4pt}
\begin{center}
\begin{tabular}{lllllll}
\hline
Class &$c$& $A$&$B$&$C$&$D$ \\
\hline
$\mathbf{C}(c)$&$\leq 3$&$B$&$\bigwedge_k\shift k^c$\\
$\mathbf{S}$&$2$&$B\ltimes W$&$k$\\
$\mathbf{T}$&$3$&$B\ltimes W$&$C\ltimes\shift(C/C_{\geq 2})$&$\bigwedge_k\shift k^2$\\
$\mathbf{B}$&$3$&$B\ltimes W$&$C\ltimes\shift C_+$&$\bigwedge_k\shift k^2$\\
$\mathbf{G}(r)$&$3$&$B\ltimes W$&$C\ltimes\Hom[k]{C}{\shift^3k}$&$k\ltimes\shift k^r$\\
$\mathbf{H}(p,q)$&$3$&$B\ltimes W$&$\Otimes[k]CD$&$k\ltimes(\shift k^p\oplus\shift^2k^q)$&$k\ltimes\shift k$\\
\hline
\end{tabular}
\end{center}\vspace{1mm}
where $W$ is a finitely generated positively graded $k$-vector space with $B_+W=0$ and $\ltimes$ designates the trivial extension from~\ref{para20201021a}.
The ring $R$ is in class $\mathbf{S}$ (that is, $A$ is of the form $k\ltimes W$) if and only if $R$ is Golod; see~\cite{golod:hslr}. If $R$ is in class $\mathbf{C}(c)$, then $R$ is a complete intersection.
\end{para}

The next step in the proof of Theorem~\ref{thm20200718b} is to connect the Poincar\'e and Bass series of $R$ to analogous series for $A$:
\begin{align*}
I_R(t)&=t^e\cdot I_A(t)&
P_k^R(t)&=(1+t)^{e}\cdot P_k^A(t)
\end{align*}
where $I_A(t)$ and $P_k^A(t)$ are the Bass series and the Poncar\'{e} series for $A$ which are defined in the DG setting as in~\eqref{eq20200902a}. These equalities are based on work in~\cite{MR485906, avramov:bsolrhoffd}.

The third step in the proof of Theorem~\ref{thm20200718b} is to analyze the Poincar\'e and Bass series of $A$ to draw the following conclusions about the corresponding series for $R$;
the proof then concludes from an analysis of the coefficients in the displayed series.

\begin{thm}[\protect{\cite[Theorem 2.1]{avramov:cslrec3}}]\label{thm20200718x}
Use the notation from~\ref{para20200816z}. Assume that $c\leq 3$ and set $l:=\rank_k A_1-1$, $n:=\rank_k A_3$, $p:=\rank_k (A_1)^2$, $q:=\rank_k (A_1\cdot A_2)$, and $r:=\rank_k(\delta_2)$, where $\delta_2\colon A_2\to \Hom[k]{A_1}{A_3}$ is defined by $\delta_2(a_2)(a_1):=a_2a_1$ for all $a_1\in A_1$ and $a_2\in A_2$.
Then the following equalities hold for the Poincar\'{e} series and Bass series of $R$:
$$
P_k^R(t)=\frac{(1+t)^{e-1}}{g(t)}\qquad \text{and}\qquad I_R^R(t)=t^d\cdot \frac{f(t)}{g(t)}
$$
where $f(t), g(t)\in \mathbb{Z}[t]$ are described as follows, where $p+q\geq 1$:\vspace{4pt}
\begin{center}
\begin{tabular}{lllll}
\hline
Class & $g(t)$&&$f(t)$& \\
\hline
$\mathbf{C}(c)$&$(1-t)^c(1+t)^{c-1}$&&$(1-t)^c(1+t)^{c-1}$&\\
$\mathbf{S}$&$1-t-lt^2$&&$1+t-t^2$&\\
$\mathbf{T}$&$1-t-lt^2-(n-3)t^3-t^5$&&$n+lt-2t^2-t^3+t^4$\\
$\mathbf{B}$&$1-t-lt^2-(n-1)t^3+t^4$&&$n+(l-2)t-t^2+t^4$\\
$\mathbf{G}(r)$&$1-t-lt^2-nt^3+t^4$&&$n+(l-r)t-(r-1)t^2-t^3+t^4$\\
$\mathbf{H}(0,0)$&$1-t-lt^2-nt^3$&&$n+lt+t^2-t^3$\\
$\mathbf{H}(p,q)$&$1-t-lt^2-(n-p)t^3+qt^4$&&$n+(l-q)t-pt^2-t^3+t^4$\\
\hline
\end{tabular}
\end{center}\vspace{1mm}
\end{thm}

We end this section with the discussion of some properties of the class $\mathbf{G}(r)$ including recent counterexamples to a conjecture of Avramov~\cite{avramov:cslrec3}.
Consider the notation from~\ref{para20200816z}. Let $n$ denote the minimal number of generators for $I$, and assume that $c=3$. In case that $R$ is a Gorenstein ring which is not complete intersection, it is known from work of J. Watanabe~\cite{MR319985} that $n\geq 5$ and $n$ is odd. Also, in this case, $R$ belongs to the class $\mathbf{G}(2i+1)$ for some $i\geq 2$ by~\cite{avramov:cslrec3}. In particular, $R$ belongs to the class $\mathbf{G}(n)$.

Conversely, Avramov {\it op.\ cit.} conjectured that
if $R$ is in the class $\mathbf{G}(r)$ with $r\geq 2$, then $R$ is Gorenstein and therefore, the classes $\mathbf{G}(3)$ and $\mathbf{G}(2i)$ for all $i\geq 1$ are empty. Christensen, Veliche, and Weyman~\cite{MR3160716, MR4038053} gave counterexamples to this conjecture. More precisely, it is shown in the latter paper that if $S$ is the power series algebra in three variables over a field, then for every $r\geq 3$ there is an ideal $I$ of $S$ with $\type(S/I)=2$ such that $S/I$ belongs to $\mathbf{G}(r)$. For counterexamples to Avramov's conjecture of arbitray type, see VandeBogert~\cite{kellerV}.

\section{Friendliness And Persistence of local rings}\label{sec3}

In this section, let $(R,\fm,k)$ be a local ring.

\subsection*{Vanishing of Ext and Tor, and finiteness of homological dimensions}\

Let $M,N$ be finitely generated $R$-modules. Following Avramov, Iyengar, Nasseh, and Sather-Wagstaff~\cite{avramov:phcnr}, $R$ is called \emph{Tor-friendly} if $\Tor [R]{i} MN=0$ for all $i\gg 0$ implies that $\pd_R M<\infty$ or $\pd_R N<\infty$. We say that $R$ is \emph{Tor-persistent} if $\Tor[R]{i}MM=0$ for all $i\gg 0$ implies that $\pd_R M<\infty$. The ring $R$ is \emph{Ext-friendly} if $\Ext[R]{i}MN=0$ for all $i\gg 0$ implies that $\pd_R M<\infty$ or $\id_R N<\infty$, where $\id$ is the injective dimension. Finally, $R$ is \emph{Ext-persistent} if $\Ext[R]{i}MM=0$ for all $i\gg 0$ implies $\pd_R M<\infty$ or $\id_R M<\infty$.

Friendliness and persistence have been studied in numerous works; see for instance~\cite{avramov:svcci, avramov:edcrcvct, avramov:phcnr, huneke:voeatoscmlr, huneke:vtci, jorgensen:gabf, jorgensen:fpdve, jorgensen:nccgr, jorgensen:itrcm, nasseh:vetfp, nasseh:ahplrdmi, nasseh:lrqdmi, nasseh:oeire, MR1974627, sega:stfcar}. The main motivation for this section is the following result in which the proofs of parts (a), (b), (c), (e), and (f) use DG algebra techniques.

\begin{thm}[\protect{\cite[Theorem 5.1, Lemmas 5.7 and 5.9]{avramov:phcnr}}]\label{thm20200803c}
Assume there exist a local homomorphism $R\to R'$ of finite flat dimension and a
deformation $R'\twoheadleftarrow Q$, i.e., a local surjection with kernel generated by a $Q$-regular sequence, where $Q$ satisfies at least one of the conditions
  \begin{enumerate}[\quad\rm(a)]
 \item
$\edim Q-\depth Q\leq 3$.
 \item
$Q$ is Gorenstein and $\edim Q-\depth Q=4$.
 \item
$Q$ is Cohen-Macaulay, almost complete intersection, $\edim Q-\depth Q=4$, and $\frac12\in Q$.
 \item
$Q$ is complete intersection.
 \item
$Q$ is one link from a complete intersection.
 \item
$Q$ is two links from a complete intersection and is Gorenstein.
 \item
$Q$ is Golod.
  \item
$Q$ is Cohen-Macaulay and $\mult Q\leq 7$.
  \end{enumerate}
Then $R$ is Tor- and Ext-persistent.
Moreover, $Q$ can be chosen to be complete, with algebraically closed residue field, and with no embedded deformation; in this case, $Q$ is Tor-friendly.
\end{thm}

One of the most important motivations for working on friendliness and persistence is the following conjecture that is known as the \emph{Auslander-Reiten Conjecture}~\cite{auslander:gvnc}. This conjecture stems from work of Nakayama~\cite{MR104718} and Tachikawa~\cite{MR0349740} on the representation theory of Artin algebras.

\begin{conj}[\protect{\cite[p.\ 70]{auslander:gvnc}}]\label{conj20200803a}
Let $M$ be a finitely generated $R$-module that satisfies the condition $\Ext[R]{i}{M}{M\oplus R}=0$ for all $i>0$. Then $M$ is a free $R$-module.
\end{conj}

\begin{para}\label{para20200803a}
It is straightforward to show that if $R$ is Ext-persistent, then it satisfies the Auslander-Reiten Conjecture~\ref{conj20200803a}.

By~\cite[Proposition 6.5]{avramov:phcnr}, Tor-friendliness implies Ext-friendliness. (In the context of complexes, these two notions are equivalent; see~\cite[Propositions 3.2 and~6.5]{avramov:phcnr}.) The question of whether all rings are Tor-persistent is open. However, examples of rings that are not Ext-persistent (hence, not Ext-friendly nor Tor-friendly) are straightforward to construct: for instance, $\Otimes[k]{(k[x,y]/(x,y)^2)}{(k[u,v]/(u,v)^2)}$.
\end{para}

Next, we describe some DG methods  from~\cite{avramov:phcnr, avramov:htecdga} used to prove Theorem~\ref{thm20200803c}.

\subsection*{Perfect DG modules, trivial extensions, and DG syzygies}\

In order to apply DG techniques in the above setting, the first tool we need is the following DG analogue of finitely generated module of finite projective dimension.

\begin{para}\label{para20200913a}
Assume that $(B, \fm_B)$ is a local DG algebra. A homologically finite DG $B$-module $M$ is called \emph{perfect} is it satisfies one of the following equivalent conditions (see~\cite{avramov:htecdga} or~\cite{test}):
\begin{enumerate}[\rm(i)]
\item
$M$ is quasiisomorphic to a DG $B$-module $F$ such that the underlying graded $B^{\natural}$-module $F^{\natural}$ has a finite basis.
\item
For all homologically bounded DG $B$-modules $N$, one has $\Tor[B]{i}MN=0$ for all $i\gg 0$.
\item
$\Tor[B]{i}M{B/\fm_B}=0$ for all $i\gg 0$.
\end{enumerate}
\end{para}

The approach described below to understanding friendliness and persistence is motivated in parts by work of Nasseh and Yoshino~\cite{nasseh:oeire} who prove that the trivial extension $R\ltimes k$ is Tor-friendly. See~\ref{para20201021a} for the definition of trivial extensions. This result is generalized to the DG setting as follows.

\begin{thm}[\protect{\cite[Theorem 4.1]{avramov:htecdga}}]\label{thm20200803a}
Let $A$ be a DG algebra that is quasiisomorphic to $B\ltimes W$, where $B$ is a homologically bounded local DG algebra, and  $W$ is a homologically bounded DG $k$-module with $\HH(W)\neq 0$.
If $M, N$ are homologically finite DG $A$-modules with $\Tor[A]{i}MN=0$ for all $i\gg 0$, then $M$ or $N$ is perfect.
\end{thm}

The proof of Theorem~\ref{thm20200803a} is similar to that of~\cite[Theorem 3.1]{nasseh:oeire}. In order to translate {\it loc.\ cit.} to the DG setting, a DG version of the important notion of a syzygy was needed. This is the DG module $N$ in the following result which we expect to be useful for other applications.

\begin{prop}[\protect{\cite[Proposition 4.2]{avramov:htecdga}}]\label{prop20200803a}
Let $(A,A_+)$ be a local DG $R$-algebra. Let $M$ be a homologically
finite DG $A$-module. Then there exists a short exact sequence
$$
0\to N\xra{\alpha} F\to \widetilde{M}\to 0
$$
of morphisms of DG $A$-modules such that
\begin{enumerate}[\rm(1)]
\item
$M\simeq \widetilde{M}$;
\item
the underlying graded $A^{\natural}$-module $F^{\natural}$ has a finite basis; and
\item
$\im(\alpha)\subseteq A_+\cdot F$.
\end{enumerate}
\end{prop}

%\begin{prop}\label{thm160527e}
%Let $A=B\ltimes W$, where $(B,B_+,k)$ is a local DG $R$-algebra and $W$ is a DG $k$-module. For each pair $(M,N)$ of DG $B$-modules there is an isomorphism:
%\[
%\Tor[A]{}MN \cong \Tor[B]{}MN \oplus \left(\Tor[B]{}Mk \otimes_{k} \Sigma \HH(W) \otimes_{k}\Tor[A]{}kN\right).
%\]
%\end{prop}
\subsection*{Friendliness and persistence}\

An important consequence of Theorem~\ref{thm20200803a} is the following result that is a bridge between Ext vanishing over $R$ and its corresponding DG algebra.

\begin{thm}[\protect{\cite[Theorem 6.3]{avramov:htecdga}}]\label{thm20200803b}
Assume there exists a minimal Cohen presentation $\widehat{R}\cong P/I$ such that
the minimal free resolution of $\widehat{R}$ over $P$ has the structure of a DG algebra and
the $k$-algebra $A=\Tor[P]{}{\widehat{R}}k$ is isomorphic to the trivial extension $B\ltimes W$ of a
graded $k$-algebra $B$ by a graded $B$-module $W\neq 0$ with $B_{\ges 1}\cdot W=0$.
Then $R$ is Tor-friendly.
\end{thm}

The proof of this result, which we outline next, relies on Avramov's machine~\ref{para20200803g} whence we also take our notation. To prove Theorem~\ref{thm20200803b}, one transfers Tor-vanishing over $R$ to Tor-vanishing over the Koszul complex $K=K^{R}(\mathbf{x})$ by base change. Then using the quasiisomorphisms~\eqref{eq201208c}, one transfers Tor-vanishing over $K$ to Tor-vanishing over $A$. Since the property of being perfect transfers from  $A$ to $K$, then to $R$, the DG result Theorem~\ref{thm20200803a} gives us the desired conclusion.

Next we sketch the proof of Theorem~\ref{thm20200803c}.
Using standard base-change techniques, one can assume without loss of generality that $R=R'$ and hence, $R$ and $Q$ have a common residue field $k$. Furthermore, we can assume that $Q$ is complete, $k$ is algebraically closed, and $Q$ does not admit embedded deformation; see~\cite[Lemma~5.7]{avramov:phcnr}. It suffices by~\cite[Theorems 2.2 and 6.3]{avramov:phcnr} and~\ref{para20200803a} to show that $Q$ is Tor-friendly.
Let $\widehat{Q}\cong P/J$ be a minimal Cohen presentation, and let $F$ be a minimal free resolution of $\widehat{Q}$ over $P$.
If $Q$ satisfies one of the conditions (a)--(g) in Theorem~\ref{thm20200803c}, then $F$ admits
a DG-algebra structure as we mentioned in~\ref{para20200718a}.
For some of these cases, the Tor algebra $\Tor[P]{}{\widehat{Q}}k$ satisfies the assumptions of Theorem~\ref{thm20200803b}. Hence, $Q$ is Tor-friendly in those cases by Theorem~\ref{thm20200803b}. In the remaining cases other methods are used to conclude that $Q$ is Tor-friendly.

Geller\cite{Hugh} and Morra~\cite{Todd} are working to apply Theorem~\ref{thm20200803b} to other rings.

\section{Bass series of local ring homomorphisms of finite flat dimension}\label{sec20200816b}

In this section, let $\varphi\colon (R,\fm, k)\to (S, \fn, \ell)$ be a local ring homomorphism.
%We start with the notion of a DG fiber which was introduced in~\cite{MR749041}.

%Avramov, Foxby, and Lescot prove the following results which are about rings and modules, however, their proofs are completely built on a heavy use of the DG techniques.

\subsection*{Relations among Bass series}\

Assume in this paragraph that $\varphi$ is flat. Then many properties of $S$ are controlled by the corresponding properties for $R$ and the closed fibre\footnote{or ``fiber,'' depending on your preference} $S/\fm S$. For instance, $S$ is Gorenstein if and only if $R$ and $S/\fm S$ are both Gorenstein. More generally, the Bass series of $S$ is related to the Bass series for $R$ and $S/\fm S$ by the formula
\begin{equation}\label{eq20201028a}
I_{S}(t)=I_{R}(t)I_{S/\fm S}(t).
\end{equation}
In particular,
for each $i\in \mathbb{Z}$, we have $\mu_R^{i+\depth R}\leq \mu_S^{i+\depth S}$. If $S/\fm S$ is Gorenstein, then Grothendieck says that $\varphi$ is Gorenstein~\cite[7.3.1--7.3.2]{grothendieck:ega4-2}.

When $\varphi$ is not flat, the properties in the previous paragraph can fail, e.g., for the natural surjection $R\to k$ when $R$ is not regular, i.e., when $\pd_R k$ is not finite. However, Avramov, Foxby, and Lescot~\cite{avramov:glh, avramov:lgh, avramov:bsolrhoffd} recognized that the full strength of flatness is not needed:

\begin{thm}[\protect{\cite[Theorems A, B, C]{avramov:bsolrhoffd}}]\label{thm20200824a}
Assume that $\varphi$ is of finite flat dimension, i.e., the $R$-module $S$ has a bounded resolution by flat modules. For instance, this holds if $S=R/I$, where $I$ is an ideal of $R$ with finite projective dimension.
\begin{enumerate}[\rm(a)]
\item
There is a formal Laurent series $I_{\varphi}(t)$ with non-negative integer coefficients such that
\begin{equation}\label{eq20200706a}
I_{S}(t)=I_{R}(t)I_{\varphi}(t).
\end{equation}
\item
For each $i\in \mathbb{Z}$, the following inequality holds:
\begin{equation*}%\label{eq20200824a}
\mu_R^{i+\depth R}\leq \mu_S^{i+\depth S}.
\end{equation*}
\item
Assume further that the closed fibre $S/\fm S$ is artinian, and either $\varphi$ is not flat or $S/\fm S$ is not a field. Then the following coefficient-wise inequality holds:
\begin{equation}\label{eq20200824b}
I_{S}(t)\preccurlyeq I_{R}(t)\frac{-(1+t)+\sum_{i=0}^{\fd_R(S)}\len_S\left(\Tor[R]{i}kS\right)t^{-i}}{1+t-\sum_{i=0}^{\fd_R(S)}\len_S\left(\Tor[R]{i}kS\right)t^{i+1}}.
\end{equation}
Equality in~\eqref{eq20200824b} holds if and only if $\varphi$ is Golod; see~\ref{para20200809s}.
\end{enumerate}
\end{thm}

\begin{para}
Here is some perspective on Theorem~\ref{thm20200824a}(c). If the closed fibre $S/\fm S$ is artinian, then the following coefficient-wise inequality holds:
\begin{equation}\label{eq20200824r}
P^{S}_{\ell}(t)\preccurlyeq \frac{P^{R}_k(t)}{1+t-\sum_{i=0}^{\fd_R(S)}\len_S\left(\Tor[R]{i}kS\right)t^{i+1}}.
\end{equation}
The ring homomorphism $\varphi$ is called a \emph{standard Golod homomorphism} if equality holds in~\eqref{eq20200824r}.

Assume either $\varphi$ is not flat or $S/\fm S$ is not a field. Then $\varphi$ is a Golod homomorphism if and only if it is a standard Golod homomorphism; see Avramov~\cite{MR846439}. Hence, in the finite flat dimension setting, Theorem~\ref{thm20200824a}(c) says that equality in~\eqref{eq20200824b} holds if and only if equality in~\eqref{eq20200824r} holds if and only $\varphi$ is Golod.
\end{para}

The proof of Theorem~\ref{thm20200824a} uses the DG fibre  introduced by Avramov~\cite{MR749041}.

\subsection*{The DG fibre of $\varphi$}\

Assume that $\varphi$ is of finite flat dimension. Let $G\xra{\simeq} k$ and $L\xra{\simeq} S$ be DG algebra resolutions over $R$. (Note that the free modules in $L$ will not be finitely generated over $R$ in general.) The \emph{DG fibre} of $\varphi$ is defined to be the local DG algebra $$F(\varphi):=\Otimes[R]GS\simeq \Otimes[R]GL\simeq \Otimes[R]kL$$ where the quasiisomorphisms come from the balance property for $\Tor[R]{}kS$. The multiplication on $F(\varphi)$ is inherited from $G$, $S$, $k$, and $L$.
The degree $0$ homology module of $F(\varphi)$ is the closed fibre $S/\fm S$. In case that $\varphi$ is flat, $F(\varphi)\simeq S/\fm S$. 

The \emph{Bass series} of $\varphi$, denoted $I_{\varphi}(t)$, is the Bass series $I_{F(\varphi)}(t)$ of the DG algebra $F(\varphi)$, which by~\cite[Theorem A]{avramov:bsolrhoffd} is a formal Laurent series.

In the case where $\varphi$ is flat, the formulas~\eqref{eq20201028a} and~\eqref{eq20200706a} are the same. In this case, they are a particular instance of the formula
$$
I^{\Otimes[R]MS}_S(t)=I^M_R(t)I_{S/\fm S}(t)
$$
where $M$ is finitely generated over $S$; one verifies this formula using the isomorphism
$$
\Ext[S]{}{\ell}{\Otimes[R]MS}\cong \Ext[R]{}kM\otimes_k \Ext[S/\fm S]{}{\ell}{S/\fm S}
$$
In the general finite flat dimension case,
Theorem~\ref{thm20200824a}(a) follows from a similar isomorphism. The innovative point in~\cite{avramov:bsolrhoffd} that we want to emphasize here is the replacement of the usual closed fibre $S/\fm S$ by the DG fibre $F(\varphi)$.

It is worth noting that Avramov and Foxby~\cite{avramov:rhafgd} established the conclusions of Theorem~\ref{thm20200824a} for a larger class of local ring homomorphisms using relative dualizing complexes, but this work does not use DG techniques.

\subsection*{Gorenstein homomorphisms}\

As we mentioned above, if $\varphi$ is flat with Gorenstein closed fibre, then $S$ is Gorenstein if and only if $R$ is Gorenstein. In case $\varphi$ has finite flat dimension, one should not expect Gorensteinness of the closed fibre to guarantee the same conclusion. In part to remedy this, Avramov and Foxby~\cite{avramov:glh, avramov:lgh} extend Grothendieck's aforementioned notion of a Gorenstein homomorphism:

The local ring homomorphism $\varphi$ is called \emph{Gorenstein} if there is an integer $a$ such that for all $i$ we have $\mu_R^{i}=\mu_S^{i+a}$. In particular, if $\varphi$ is Gorenstein, then $S$ is Gorenstein if and only if $R$ is Gorenstein. If $\varphi$ has finite flat dimension, Gorensteinness of $\varphi$ is equivalent to having the equality $\mu_R^{i}=\mu_S^{i+\depth S-\depth R}$ for all $i$ by Theorem~\ref{thm20200824a}(a).

In case that $\varphi$ is flat, Gorensteinness of $\varphi$ is equivalent to the Gorensteinness of the closed fibre $S/\fm S$; see~\cite[(4.2) Proposition]{avramov:lgh}. Hence, this notion of Gorenstein homomorphisms is a generalization of Grothendieck's Gorenstein homomorphisms.

The result {\it op.\ cit.} can be extended to the following characterization of Gorenstein homomorphisms in terms of their DG fibres.

\begin{thm}[\protect{\cite[(4.4) Theorem]{avramov:lgh}}]\label{thm20200824t}
Assume that $\varphi$ has finite flat dimension. Then $\varphi$ is Gorenstein if and only if the DG fibre $F(\varphi)$ is a Gorenstein DG algebra (that is, $I_{\varphi}(t)=t^d$ for some integer $d$).
\end{thm}

As one might imagine, given the usefulness of the Gorenstein property for local rings, Gorenstein DG algebras have been investigated separately; see Frankild, Iyengar, and J\o rgensen~\cite{frankild:ddgmgdga, frankild:gdga}.

\section{Ascent property of pd-test modules}\label{sec20200816a}

In this section, let $\varphi\colon (R, \fm, k)\to (S, \fn, \ell)$ be a flat local ring homomorphism.

\subsection*{Pd-test modules}\

A useful, classical result states that the residue field $k$ has the ability to test for finite projective dimension: a finitely generated $R$-module $N$ has finite projective dimension if and only if $\Tor[R]{i}kN=0$ for $i\gg 0$. According to the following definition, which was coined by O. Celikbas, Dao, and Takahashi~\cite{CDtest}, this says that $k$ is a pd-test $R$-module.

A finitely generated $R$-module $M$ is called a \emph{pd-test module} if for every finitely generated $R$-module $N$ with $\Tor[R]{i}MN=0$ for $i\gg 0$ we have $\pd_R N<\infty$.

It is natural to ask how the pd-test property for a finitely generated $R$-module $M$ behaves under completion. This is related to the well-known fact that $R$ is regular if and only if $\widehat{R}$ is regular. It is straightforward to show that if $\widehat{M}$ is pd-test over $\widehat{R}$, then $M$ is pd-test over $R$. That is, the pd-test property descends from the completion. The question of ascent is more subtle. It was posed in~\cite{CDtest} and answered by O. Celikbas and Sather-Wagstaff~\cite{celikbas:tgp} using derived category techniques. The following more general ascent result is proved by Sather-Wagstaff~\cite{test}.

\begin{thm}[\protect{\cite[Theorem 4.8]{test}}]\label{thm20200820a}
Assume that the closed fibre $S/\fm S$ of $\varphi$ is regular and the induced field extension $k\to \ell$ is algebraic. If a finitely generated $R$-module $M$ is pd-test over $R$, then $S\otimes_RM$ is a pd-test module over $S$.
\end{thm}

Theorem~\ref{thm20200820a} is proved using the following DG techniques.

\subsection*{Pd-test DG modules}\

A homologically finite DG module $M$ over a local DG algebra $B$ is a \emph{pd-test DG module} if every homologically finite DG $B$-module $N$ with $\Tor[B]{i}MN=0$ for all $i\gg 0$ is perfect.

The following result is a special case of a DG version of Theorem~\ref{thm20200820a}. It plays an essential role in the proof of Theorem~\ref{thm20200820a}.

\begin{thm}[\protect{\cite[Theorem 4.6]{test}}]\label{thm20201013x}
Let $A$ be a finite-dimensional DG $k$-algebra with $A_0=k$ and $\HH_0(A)\neq 0$. Let $k\to \ell$ be an algebraic field extension, and set $B=\ell\otimes_k A$. If $M$ is pd-test over $A$, then $\Otimes[A]BM$ is pd-test over $B$.
\end{thm}

Before applying Theorem~\ref{thm20201013x}, we sketch its proof.
Assume that $N$ is a homologically finite DG $B$-module such that $\Tor[B]{i}{\Otimes[A]BM}{N}=0$ for all $i\gg 0$.
In case that $k\to \ell$ is a finite field extension, the assertion follows from a standard argument using~\ref{para20200913a}.
Now consider the general case, where $k\to \ell$ is algebraic. By truncating an appropriate resolution of $N$ over $B$ one can assume that $N$ is finite-dimensional over $\ell$. It then follows that the differential and scalar multiplication on $N$
are represented by matrices consisting of finitely many elements of $\ell$. Adjoining these algebraic elements to $k$, one obtains an intermediate field extension $k\to k'\to \ell$ such that $k\to k'$ is  finite. By construction of $k'$, with $A'=k'\otimes_kA$, there is a bounded DG $A'$-module $L$ such that $N\cong B\otimes_{A'}L$. At this point, the assumption of $\Tor[B]{i}{\Otimes[A]BM}{N}=0$ for all $i\gg 0$ implies that $\Tor[A']{i}{\Otimes[A]{A'}{M}}{L}=0$ for all $i\gg 0$. Since $k\to k'$ is finite, it follows that $L$ is perfect over $A'$, so $N\cong B\otimes_{A'}L$ is perfect over $B$.

\subsection*{Outline of the proof of Theorem~\ref{thm20200820a}}\

Assume that $M$ is a pd-test module over $R$. We need to show that $S\otimes_R M$ is a pd-test module over $S$. Assume that $\Tor[S]{i}{S\otimes_R M}{N}=0$ for $i\gg 0$, where $N$ is a finitely generated $S$-module.
Standard techniques reduce to the case where $R$ and $S$ are complete with $S/\fm S =\ell$. Using the notation from~\ref{para20200803g} and applying~\cite[(1.6) Theorem]{avramov:solh} we have a minimal Cohen presentation $P'\xra{\tau'} S$ and a commutative diagram of local ring homomorphisms
\begin{equation*}%\label{eq20200901b}
\begin{split}
\xymatrix{
P\ar[r]^{\alpha}\ar[d]_{\tau}&P'\ar[d]^{\tau'}\\
R\ar[r]^{\varphi}&S
}
\end{split}
\end{equation*}
such that $\alpha$ is flat, $\tau'$ is surjective, $P'/\fp P'\cong \ell$, and $S\cong R\otimes_P P'$. The last isomorphism implies that $F':=F\otimes_PP'\xra{\simeq} S$ is a DG algebra resolution of $S$ over $P'$. Note that $\varphi(\mathbf{x})$ minimally generates $\fn$. Following the process of~\ref{para20200803g} for the ring $S$, we get the next commutative diagram of morphisms of DG algebras
 \begin{equation*}%\label{eq20200901a}
\begin{split}
\xymatrix{
R\ar[r]\ar[d]_{\varphi}&K^R\ar[d]&K^P\otimes_P R\ar[l]_{\cong}\ar[d]&K^P\otimes_P F\ar[l]_{\simeq}\ar[d]\ar[r]^{\simeq}& k\otimes_P F=A\ar[d]\\
S\ar[r]&K^S&K^{P'}\otimes_{P'} S\ar[l]_{\cong}&K^{P'}\otimes_{P'} F'\ar[l]_{\simeq}\ar[r]^{\simeq}& \ell\otimes_{P'} F'
}
\end{split}
\end{equation*}
in which $K^R=K^R(\mathbf{x})$, $K^S=K^S(\varphi(\mathbf{x}))$, $K^{P}=K^{P}(\mathbf{y})$, and $K^{P'}=K^{P'}(\alpha(\mathbf{y}))$.
Note that the DG algebra $\ell\otimes_{P'} F'$ is isomorphic to $\ell\otimes_k A$. Now, the pd-test problem between $R$ and $S$ can be translated through the rows of this diagram to a DG pd-test problem between $A$ and $\ell\otimes_k A$. At this point the assertion follows from Theorem~\ref{thm20201013x}.

In case that $\ell=k(x)$ is a transcendental extension of $k$, the same conclusion as in the statement of Theorem~\ref{thm20200820a} holds by a result of Tavanfar~\cite{Ehsan}.

\section{A conjecture of Vasconcelos on the conormal module}

Throughout this section, let $I$ be an ideal of a ring $R$, and set $S=R/I$. 

Ferrand~\cite{MR219546} and Vasconcelos~\cite{MR213345} show that properties of the ring $S$ are often reflected in the properties of the \emph{conormal module} $I/I^2$ over $S$. This section focuses on the following conjecture of Vasconcelos~\cite{MR508082}.

\begin{conj}[\protect{\cite[($C_1$)]{MR508082}}]\label{conj20201012a}
If $\pd_R S$ and $\pd_{S}I/I^2$ are finite, then $I$ is locally generated by a regular sequence.
\end{conj}

This conjecture was settled in the affirmative for some special cases by Vasconcelos~\cite{MR814190}, Gulliksen and Levin~\cite{MR0262227}, and Herzog~\cite{herzogv}. The following major progress on this conjecture was made by Avramov and Herzog~\cite{MR1269426} using Andr\'{e}-Quillen homology and DG homological methods.

\begin{thm}[\protect{\cite[Theorem 3]{MR1269426}}]\label{thm20201012c}
Let $k$ be a field of characteristic $0$, and assume $R$ is a positively graded polynomial ring over $k$ and $I$ is homogeneous. Then the following are equivalent:
\begin{enumerate}[\rm(i)]
\item
$S$ is complete intersection;
\item
$I/I^2$ is a free $S$-module;
\item
$\pd_{S}I/I^2<\infty$.
\end{enumerate}
\end{thm}

In a recent paper, Briggs~\cite{briggs2020vasconcelos} establishes Conjecture~\ref{conj20201012a} in its full generality.

\begin{thm}[\protect{\cite[Theorem A]{briggs2020vasconcelos}}]\label{thm20201012d}
Conjecture~\ref{conj20201012a} holds in general.
\end{thm}

\begin{para}\label{para20201031s}
Briggs' proof for Theorem~\ref{thm20201012d} relies on methods pioneered by Avramov and Halperin~\cite{MR749041, avramov:tlg} on homotopy Lie algebras $\pi^*(\varphi)$ arising from DG constructions.

Assume without loss of generality that $(R,\fm,k)$ and $(S,\fn,k)$ are local. Let $\varphi\colon R\to S$ be the natural surjection. Fix a \emph{minimal model} for $\varphi$ which is a factorization $R\to A\xra{\simeq}S$, where $(A, \fm_A)$ is a local DG $R$-algebra such that:
\begin{enumerate}[\rm(a)]
\item
The underlying algebra $A^{\natural}=R[X_1,X_2,\ldots]$ is the free graded commutative $R$-algebra, where each $X_i$ is a set of variables of degree $i$; and
\item
$\partial(\fm_A)\subseteq \fm+\fm_A^2$.
\end{enumerate}
The DG algebra $A$ is also denoted $R\langle X\rangle$.

A graded basis for each $\pi^i(\varphi)$ is dual to $X_i$, and each element $z\in\pi^2(\varphi)$ corresponds to a derivation $\theta_z\colon A\to \fm_A$ of degree $-2$ as is described in~\cite{avramov:ifr, briggs2020vasconcelos}. Let $\overline{\theta_z}\colon A\to \fn$ be the composition of $\theta_z$ and the surjective quasiisomorphism $\fm_A\to\fn$. Under the assumptions of Conjecture~\ref{conj20201012a}, one can find a certain factorization of $\overline{\theta_z}$ which implies that $z$ is radical in $\pi^2(\varphi)$; see~\cite[proof of Lemma 2.6 and Theorem~2.7]{briggs2020vasconcelos}. Now~\cite[Theorem C]{avramov:tlg} implies that $\varphi$ is complete intersection, as desired.
\end{para}

\section{A conjecture of Vasconcelos on semidualizing modules}\label{sec4}
In this section, $(R,\fm,k)$ is a local ring.

Here we discuss a class of modules that are particularly well-suited for creating dualities. They were originally introduced by Foxby~\cite{foxby:gmarm} who called them \emph{PG modules of rank $1$}. They are useful, e.g., for understanding Gorenstein dimensions, in particular, Avramov and Foxby's composition question for local ring homomorphisms of finite G-dimension~\cite{avramov:rhafgd, sather:cidfc}.

\subsection*{Semidualizing modules}\

A finitely generated $R$-module $C$ is called \emph{semidualizing}
if the homothety morphism
$\chi^R_C\colon R\to\Hom[R]{C}{C}$ is an isomorphism and $\Ext[R]{i}CC=0$ for all $i\geq 1$. A semidualizing module of finite injective dimension is called a \emph{dualizing module}.
Let $\s_0(R)$ be the
set of  isomorphism classes of semidualizing $R$-modules.

This section is centered on the following conjecture posed by Vasconcelos~\cite{vasconcelos:dtmc}.

\begin{conj}[\protect{\cite[p. 97]{vasconcelos:dtmc}}]\label{vasc}
If $R$ is Cohen-Macaulay, then $\s_0(R)$ is finite.
\end{conj}

Note that if $R$ is Ext-persistent, then $R$ satisfies this conjecture. Moreover, in this case, the only semidualizing $R$-modules are the free module of rank $1$ and a dualizing module, if one exists.

Christensen and Sather-Wagstaff~\cite{christensen:cmafsdm} answered Conjecture~\ref{vasc} in the case where $R$ contains a field. Their proof reduces to the case of a finite-dimensional algebra, then implicitly uses the following technology from geometric representation theory.

\begin{para}\label{para20201008d}
Assume that $R$ is a finite-dimensional $k$-algebra, where $k$ is algebraically closed. The $R$-modules of a fixed length $r$ are parametrized
by an algebraic variety $\od_r^R$. One can define an action of the general linear group
$\ggl_r^k$ on $\od_r^R$. The isomorphism class of an $R$-module $M$
is  the orbit $\ggl_r^k\cdot M$,
and the tangent space $\ta^{\ggl_r^k \cdot M}_M$
to the orbit $\ggl_r^k \cdot M$ at $M$
is identified with
a subspace of the tangent space $\ta^{\od_r^R}_M$.
A result of Voigt~\cite{voigt:idteag}
(see also Brion~\cite{brion:rq} or Gabriel~\cite{gabriel:frto})
provides an isomorphism $\Ext 1MM\cong \ta^{\od_r^R}_M/\ta^{\ggl_r^k \cdot M}_M$.
As in work of Happel~\cite{happel:sm}, it follows that if $\Ext 1MM=0$ (e.g., if $M$ is a semidualizing module), then the orbit $\ggl_r^k \cdot M$ is open in $\od_r^R$. Since $\od_r^R$ is quasi-compact, it can contain only finitely many open orbits, hence, $\s_0(R)$ is finite.
\end{para}

Using a modification of these ideas, Nasseh and Sather-Wagstaff~\cite{nasseh:gart} establish Conjecture~\ref{vasc} in total generality with no Cohen-Macaulay hypothesis.

\begin{thm}[\protect{\cite[Theorem A]{nasseh:gart}}]\label{main'}
For the local ring $R$, the set $\s_0(R)$ is finite.
\end{thm}

\subsection*{A DG version of Voigt's theorem and the proof of Theorem~\ref{main'}}\

To prove Theorem~\ref{main'}, we work with the following DG version of semidulazing modules due to Christensen and Sather-Wagtaff~\cite{christensen:dvke}.

Let $A$ be a homologically degreewise noetherian DG $R$-algebra. A homologically finite DG $A$-module $C$ is \emph{semidualizing}
if the homothety morphism
$\chi^A_C\colon A\to\Rhom[A]{C}{C}$ is an isomorphism in the derived category $\catd(A)$.
If $A=R$, a semidualizing DG $R$-module $C$ is called a \emph{semidualizing $R$-complex}. A semidualizing $R$-complex of finite injective dimension is called a \emph{dualizing complex}.
Let $\s(A)$ denote the set of shift-isomorphism classes
of semidualizing DG $A$-modules in $\catd(A)$.

Theorem~\ref{main'} is a consequence of the following result because $\s_0(R)\subseteq \s(R)$.

\begin{thm}[\protect{\cite[4.2 and Theorem A]{nasseh:gart}}]\label{DGmain}
Consider the notation of~\ref{para20200803g}. The sets $\s(A)$ and $\s(R)$ are finite.
\end{thm}

Using Grothendieck~\cite[Proposition (0.10.3.1)]{grothendieck:ega3-1}, we can assume in Theorem~\ref{DGmain} that $R$ is complete with algebraically closed residue field. Because of Avramov's machine~\ref{para20200803g}, it suffices to show that $\s(A)$ is finite. To establish this finiteness, one uses the following  DG version of~\ref{para20201008d} above.

The set of finite-dimensional DG $A$-modules $M$ with fixed
underlying graded $k$-vector space $W$ is parametrized
by an algebraic variety $\od^A(W)$.
A product $\ggl(W)_0$ of
general linear groups acts on $\od^A(W)$
and the isomorphism class of $M$
is the orbit $\ggl(W)_0\cdot M$ under this action. See~\cite{nasseh:gart} for more details.

The DG version of Voigt's result from~\ref{para20201008d} that enables us to prove Theorem~\ref{DGmain} is the following.

\begin{thm}[\protect{\cite[Theorem B]{nasseh:gart}}]\label{thm20200810w}
Let $W$ be a finite-dimensional graded $k$-vector space.
Given an element $M\in\od^A(W)$, there is an isomorphism
$$\ta^{\od^A(W)}_M/\Taglm_M \cong \yext[A]{1}MM$$
where $\yext[A]{1}MM$ denotes the Yoneda Ext group defined as the set of equivalence
classes of short exact sequences $0\to M\to L\to M\to 0$.
\end{thm}

As in~\ref{para20201008d}, it follows from Theorem~\ref{thm20200810w} that if $\yext[A] 1MM=0$, then the orbit $\ggl(W)_0\cdot M$ is open in $\od^A(W)$. Since $\od^A(W)$ is quasi-compact, it follows that there are only finitely many open orbits in it. Thus, it remains to show that a each semidualizing DG $A$-module $C$ satisfies $\yext[A] 1CC=0$. This vanishing follows from work of Nasseh and Sather-Wagstaff~\cite{nasseh:egdgm}.

One can actually obtain a very tight connection between the sizes of $\s(A)$ and $\s(R)$ using a lifting result in~\cite{nasseh:ldgm} that generalizes results of Auslander, Ding, and Solberg~\cite{auslander:lawlom} and Yoshino~\cite{yoshino}. See Nasseh and Yoshino~\cite{MR3774891} and Ono and Yoshino~\cite{MR4152766} for more general lifting results.
Also, Altmann and Sather-Wagstaff~\cite{Hannah} utilize Avramov's machine to extend results of Gerko~\cite{gerko:sdc} from the realm of finite-dimensional algebras to arbitrary local rings.

\section{Complete intersection maps and the proxy small property}\label{sec20200908c}

In this section, let $\varphi\colon R\to S$ be a surjective ring homomorphism. 

Here, we outline results of Briggs, Iyengar, Letz, and Pollitz~\cite{BenBriggs} on questions motivated by work of Dwyer, Greenlees, and Iyengar~\cite{MR2225632} and Pollitz~\cite{MR3988642}.

A triangulated subcategory $\mathcal{X}$ of the derived category $\mathcal{D}(R)$ is called \emph{thick} if it is closed under direct summands and satisfies the following two-of-three property: for each exact triangle $L\to M\to N\to$ in $\mathcal{D}(R)$ if two of the objects are in $\mathcal{X}$, then so is the third. The thick subcategory of $\mathcal{D}(R)$ \emph{generated by} an $R$-complex $M$ is the smallest thick subcategory of $\mathcal{D}(R)$ (with respect to inclusion) that contains $M$. Note that an $R$-complex is perfect if and only if it is in the thick subcategory generated by $R$. If an $R$-complex $N$ is in the thick subcategory generated by another $R$-complex $M$, we say that $N$ is \emph{finitely built from} $M$.

A triangulated subcategory of $\mathcal{D}(R)$ is called \emph{localizing} if it is closed under
arbitrary coproducts. Note that a localizing subcategory is thick.
The localizing subcategory of $\mathcal{D}(R)$ \emph{generated by} an $R$-complex $M$ is the smallest localizing subcategory of $\mathcal{D}(R)$ that contains $M$. If an $R$-complex $N$ is in the localizing subcategory generated by another $R$-complex $M$, we say that $N$ is \emph{built from} $M$.

A \emph{small} complex $M$ over a ring $R$ is an $R$-complex such that $\Hom[\mathcal{D}(R)]M-$ commutes with arbitrary direct sums. Note that the perfect $R$-complexes are precisely the small $R$-complexes (or the small objects in $\mathcal{D}(R)$).

In~\cite{MR2200850}, an $R$-complex $M$ is \emph{proxy small} if there exists a small $R$-complex $N$ such that
$N$ is finitely built from $M$, and
$M$ is built from $N$.
Note that every small $R$-complex is proxy small. Other examples of proxy small complexes include the residue field of a local ring and modules of finite complete intersection dimension over a local ring.

Let $R$ be a local ring. The famous result of Auslander-Buchsbaum and Serre~\cite{auslander:cm, serre:sldhdaedmn} says that $R$ is regular if and only if every homologically bounded $R$-complex is small. The paper~\cite{MR2225632} contains a partial analogue of this statement for complete intersection rings: if $R$ is  complete intersection, then every homologically bounded $R$-complex is proxy small. Pollitz~\cite{MR3988642} proved the converse of this by showing that if every homologically bounded $R$-complex is proxy small, then $R$ is complete intersection. Pollitz's proof heavily uses DG methods relying on his version~\cite{Josh2} of Avramov and Buchweitz's~\cite{avramov:svcci} support varieties over Koszul complexes. Due to space restrictions here, we do not provide further details of this construction.

In the not necessarily local setting, \cite{MR2225632} includes a more general statement than the one mentioned above: if $\varphi$ is complete intersection, then proxy smallness ascends along $\varphi$, i.e., any $S$-complex that is proxy small over $R$ is proxy small over $S$.
Briggs, Iyengar, Letz, and Pollitz~\cite{BenBriggs} prove the following converse of this statement.

\begin{thm}[\protect{\cite[Theorem B]{BenBriggs}}]\label{thm20200914c}
Assume that $\varphi$ has finite projective dimension. If proxy smallness ascends along $\varphi$, then $\varphi$ is complete intersection.
\end{thm}

A consequence of this theorem~\cite[Corollary 4.1]{BenBriggs} is another proof of a fundamental result of Avramov~\cite[(5.7.1) Lemma]{avramov:lcih} used in his solution to Quillen's conjecture discussed in Section~\ref{sec20200908a} below. More precisely, if $R\xra{\varphi}S\xra{\psi}T$ are surjective local homomorphisms such that $\fd_ST<\infty$, then $\psi\circ \varphi$ is complete intersection if and only if $\varphi$ and $\psi$ are complete intersection.

The proof of Theorem~\ref{thm20200914c} reduces to the case where $(R,\fm)$ and $(S,\fn)$ are local. Set $\widetilde{S}:=R/I$, where $I$ is an ideal generated by a maximal $R$-regular sequence in $\ker \varphi\setminus \fm\ker \varphi$. The surjection $R\to S$ is the composition of the natural surjections $R\xra{\widetilde{\varphi}}\widetilde{S}\xra{\dot{\varphi}}S$. To complete the proof it suffices to show that $S$ is small over $\widetilde{S}$; indeed, then~\cite[Corollary 1.4.7]{bruns:cmr} implies $\varphi=\widetilde{\varphi}$  is complete intersection, as desired.

To show that $S$ is small over $\widetilde{S}$, let $K=K^S(\fn)$ be the Koszul complex on a minimal generating set of $\fn$, and consider the restriction $\dot{\varphi}_*\colon \mathcal{D}(S)\to \mathcal{D}(\widetilde{S})$. By~\cite[Remark 5.6]{MR2225632}, it suffices to prove that $\dot{\varphi}_*\left(K\right)$ is a small $\widetilde{S}$-complex. This smallness follows from the next lemma which uses Hochschild cohomology for DG algebras as constructed by Avramov, Iyengar, Lipman, and Nayak~\cite{MR2565548}.

\begin{lem}[\protect{\cite[Lemma 2.5]{BenBriggs}}]\label{lem20201031a}
Let $A$ be a DG $R$-algebra, and let $M$ and $N$ be DG $A$-modules. Let $\alpha$ be an element of the graded Hochschild cohomology algebra $\HH\!\HH^*(A\!\mid\! R)$. If $N$ is (finitely) built from $M$, then the mapping cone $N/\!/\alpha$ of an induced morphism $N\xra{\chi_N(\alpha)} \shift^{|\alpha|}N$ is (finitely) built from $M/\!/\alpha$.
In particular, if $M$ is proxy small then so is $M/\!/\alpha$.
\end{lem}

\section{Conjectures of Quillen on Andr\'{e}-Quillen homology}\label{sec20200908a}

In this section, let $\varphi\colon R\to S$ be a ring homomorphism.

Here, we describe Avramov's solution~\cite{avramov:lcih} to a famous conjecture of Quillen~\cite{MR0257068} and Avramov and Iyengar's significant progress~\cite{MR1977825} on a second one.

\subsection*{Quillen's conjectures}\

The $n$th \emph{Andr\'{e}-Quillen homology} of the $R$-algebra $S$ with coefficients in an $S$-module $N$ is $\DD_n(S\!\mid\! R, N)=\HH_n\left(\LD(S\!\mid\! R)\otimes_S N\right)$, where $\LD(S\!\mid\! R)$ is the cotangent complex of $\varphi$; see Andr\'{e}~\cite{MR0214644}, Iyengar~\cite{MR2355775}, and Quillen~\cite{MR0257068} for definitions and foundational properties.
%Note that by definition, $\DD_n(S\!\mid\! R, -)=0$ for $n> m$ if and only if $\fd_S \LD(S\!\mid\! R)\leq m$.

The first of Quillen's conjectures that we consider deals with locally complete intersection homomorphisms. This notion was originally defined for maps that are essentially of finite type or flat. Avramov's solution of this conjecture hinges on the following generalization of this notion.

Assume in this paragraph that $\varphi\colon R\to (S,\fn)$ is a local ring homomorphism, and let $\grave{\varphi}\colon R\to \widehat{S}$ be the composition of $\varphi$ with the natural completion map $S\to \widehat{S}$. A \emph{Cohen factorization} of $\grave{\varphi}$ is a factorization into local ring homomorphisms $R\xra{\dot{\varphi}} R'\xra{\varphi'} \widehat{S}$ such that $\dot{\varphi}$ is flat with regular closed fibre, $\varphi'$ is surjective, and $R'$ is complete. If there is a Cohen factorization $R\to R'\xra{\varphi'} \widehat{S}$ of $\grave{\varphi}$ in which $\ker \varphi'$ is generated by an $R'$-regular sequence, then $\varphi$ is called \emph{complete intersection at $\fn$}.

In general, the (not necessarily local) ring homomorphism $\varphi\colon R\to S$ is called \emph{locally complete intersection} if it is complete intersection at all prime ideals $\fq$ of $S$, i.e., for all such $\fq$, the induced local ring homomorphism $\varphi_{\fq}\colon R_{\fq\cap R}\to S_{\fq}$ is complete intersection at $\fq S_{\fq}$. Also, $\varphi$ is \emph{locally of finite flat dimension} if $\fd_{R} S_{\fq}<\infty$ for all prime ideals $\fq$ of $S$. In case that $R$ has finite Krull dimension this condition is equivalent to $\fd_RS<\infty$; see Auslander and Buchsbaum~\cite{MR96720}.

Now we can state the conjectures of Quillen~\cite{MR0257068} that we are concerned with.

\begin{conj}[\protect{\cite[(5.6) and (5.7)]{MR0257068}}]\label{conj20200909a}
Assume $\varphi$ is essentially of finite type.
\begin{enumerate}[\rm(a)]
\item
If $\varphi$ is locally of finite flat dimension and $\DD_n(S\!\mid\! R, -)=0$ for all $n\gg 0$, then it is locally complete intersection.
\item
If $\DD_n(S\!\mid\! R, -)=0$ for all $n\gg 0$, then $\DD_n(S\!\mid\! R, -)=0$ for all $n\geq 3$.
\end{enumerate}
\end{conj}

\subsection*{Avramov's solution of Conjecture~\ref{conj20200909a}(a) via DG techniques}

\begin{thm}[\protect{\cite[(1.3)]{avramov:lcih}}]\label{thm20200909a}
Conjecture~\ref{conj20200909a}(a) holds without the essentially of finite type assumption.
\end{thm}

The proof of Theorem~\ref{thm20200909a} reduces to the case where $\varphi$ is surjective and local. In this case, the proof hinges on the following spectral sequence~\cite[(4.2) Theorem]{avramov:lcih}
$$
{}^2\!E_{p,q}=\pi_{p+q}\left(\operatorname{Sym}^{\ell}_q(\shift\LD(S\!\mid\! R)\otimes_S\ell)\right)\Longrightarrow \ell\langle X\rangle_{p+q}
$$
where $\ell$ is the residue field of $S$, and the other notation including the DG algebra $\ell\langle X\rangle$ is from~\ref{para20201031s}.

Very recently Briggs and Iyengar~\cite{BenSri} improved upon Theorem~\ref{thm20200909a} with the following. The proof of this result also uses DG technology, but we do not discuss it because of space constraints.

\begin{thm}[\protect{\cite[Theorem A]{BenSri}}]\label{thm20201031f}
If $\varphi$ is locally of finite flat dimension and one has $\DD_n(S\!\mid\! R, -)=0$ for some $n\geq 1$, then $\varphi$ is locally complete intersection.
\end{thm}

\subsection*{Conjecture~\ref{conj20200909a}(b) for algebra retracts}\

Avramov and Iyengar~\cite{MR1977825} proved Conjecture~\ref{conj20200909a}(b) in the case where $S$ is an \emph{algebra retract} of $R$, that is, where there is a ring homomorphism $\psi\colon S\to R$ such that $\varphi\circ \psi= \id_S$.

\begin{thm}[\protect{\cite[Theorem I]{MR1977825}}]\label{thm20200909b}
Assume that $S$ is an algebra retract of $R$. Then the following conditions are equivalent.
\begin{enumerate}[\rm(i)]
\item
$\DD_n(S\!\mid\! R, -)=0$ for all $n\gg 0$.
\item
$\DD_n(S\!\mid\! R, -)=0$ for all $n\geq 3$.
\item
$\DD_3(S\!\mid\! R, -)=0$.
\item
$\DD_n(S\!\mid\! R, -)=0$ for some $n\geq 3$ such that $\lfloor\frac{n-1}{2}\rfloor!$ is invertible in $S$.
\end{enumerate}
\end{thm}

Conjecture~\ref{conj20200909a}(b) fails in the non-noetherian case; see Andr\'{e}~\cite{MR1418037} and Planas-Vilanova~\cite{Planas}. This conjecture is still open in general for noetherian rings.

In the proof of Theorem~\ref{thm20200909b}, the following notion plays an essential role.
A local homomorphism $\varphi\colon (R,\fm,k)\to (S,\fn,\ell)$ is \emph{almost small} if the kernel of the homomorphism
$\Tor[\varphi]{}{\overline{\varphi}}{\ell}\colon \Tor[R]{}{k}{\ell}\to \Tor[S]{}{\ell}{\ell}$
of graded algebras is generated by elements of degree $1$.

DG techniques are crucial in the proof of Theorem~\ref{thm20200909b}. Key to this is a structure theorem~\cite[4.11 Theorem]{MR1977825} for surjective almost small homomorphisms in terms of DG algebra homomorphisms. From this one concludes~\cite[5.6. Theorem]{MR1977825} that almost small homomorphisms have finite \emph{weak category}; a notion motivated by the works of F\'{e}lix and Halperin~\cite{MR664027}. As a result, information on the positivity and growth of deviations of almost small homomorphisms is revealed by~\cite[5.4. Theorem]{MR1977825}. The local version of Theorem~\ref{thm20200909b} follows from this via a characterization of complete intersection local homomorphisms having finite weak category in terms of the vanishing of the Andr\'{e}-Quillen homology with coefficients in the residue field; see~\cite[6.4. Theorem]{MR1977825}. A reduction to the local case then finishes the proof.% of Theorem~\ref{thm20200909b}.

\section{Finite Generation of Hochschild Homology Algebras}\label{sec20201018a}

Throughout this section, let $\varphi\colon R\to S$ be a ring homomorphism.

We discuss work of Avramov and Iyengar~\cite{MR1779800} on finite generation of Hochschild homology algebras. In it, they prove the converse of the Hochschild-Kostant-Rosenberg Theorem using DG methods and Andr\'{e}-Quillen homology; see~\cite{MR0077480, MR2355775, MR1217970} for definitions and facts that are used in this section.

The Hochschild homology algebra, denoted $\HH\!\HH_*(S\!\mid\! R)$, is a graded commutative algebra defined using shuffle products on the Hochschild complex. This satisfies $\HH\!\HH_0(S\!\mid\! R)=S$, and $\HH\!\HH_1(S\!\mid\! R)=\Omega^1_{S\mid R}$ is the $S$-module of K\"{a}hler differentials.
%There is a canonical homomorphism
%$$
%\omega^*_{S\mid R}\colon \wedge^*_S\Omega^1_{S\mid R}\to \HH\!\HH_*(S\mid R)
%$$
%of graded algebras that maps differential forms to Hochschild homology.
Recall that the $R$-algebra $S$ is called \emph{regular} if $\varphi$ is flat and $S\otimes_R k$ is regular for each homomorphism $R\to k$ from $R$ to a field $k$.
Hochschild, Kostant, and Rosenberg~\cite{MR142598} proved that if $R$ is a perfect field and $S$ is \emph{smooth} over $R$ (that is, $S$ is a regular $R$-algebra and essentially of finite type), then $\HH\!\HH_*(S\!\mid\! R)$ is a finitely generated $S$-algebra.
Here is the aforementioned converse.

\begin{thm}\cite[Theorem (5.3)]{MR1779800}\label{thm20200913a}
Assume that $\varphi$ is flat and essentially of finite type. If
the $S$-algebra $\HH\!\HH_*(S\!\mid\! R)$ is finitely generated, then $S$ is  smooth over $R$.
\end{thm}

This result settles a conjecture of Vigu\'{e}-Poirrier~\cite{MR1245091} who already established it in the case where $S=R[x_1,\ldots,x_n]/I$, and $R$ is a field of characteristic $0$, and $I$ is generated by a regular sequence. It was also known for positively graded $S$ such that $S_0=R$ is a field of characteristic $0$ by Dupont and Vigu\'{e}-Poirrier~\cite{MR1803231}.

The DG techniques used in the proof of Theorem~\ref{thm20200913a} are confined to the  characteristic-$0$ case. Here Avramov and Iyengar use a version of Avramov's machine~\cite[4.2]{MR1779800} which gives a DG algebra $A$ where $\HH(A)$ is the Tor algebra $\Tor[R]{}SS$.

\section*{Acknowledgments}
We are grateful to Josh Pollitz and Keller VandeBogert for helpful suggestions about this survey.

%\bibliography{+new}

\begin{thebibliography}{100}

\bibitem{Hannah}
H.~Altmann and S.~Sather-Wagstaff, \emph{Chains of semidualizing complexes}, in
  preparation.

\bibitem{MR0214644}
M.~Andr\'{e}, \emph{M\'{e}thode simpliciale en alg\`ebre homologique et
  alg\`ebre commutative}, Lecture Notes in Mathematics, Vol. 32,
  Springer-Verlag, Berlin-New York, 1967. \MR{0214644}

\bibitem{MR1418037}
\bysame, \emph{Examples of non-rigid cotangent complexes}, J. Algebra
  \textbf{186} (1996), no.~1, 32--46. \MR{1418037}

\bibitem{auslander:cm}
M.~Auslander and D.~A. Buchsbaum, \emph{Codimension and multiplicity}, Ann. of
  Math. (2) \textbf{68} (1958), 625--657. \MR{0099978 (20 \#6414)}

\bibitem{MR96720}
\bysame, \emph{Homological dimension in noetherian rings. {II}}, Trans. Amer.
  Math. Soc. \textbf{88} (1958), 194--206. \MR{96720}

\bibitem{auslander:lawlom}
M.~Auslander, S.\ Ding, and \O.\ Solberg, \emph{Liftings and weak liftings of
  modules}, J. Algebra \textbf{156} (1993), 273--397. \MR{94d:16007}

\bibitem{auslander:gvnc}
M.~Auslander and I.~Reiten, \emph{On a generalized version of the {N}akayama
  conjecture}, Proc. Amer. Math. Soc. \textbf{52} (1975), 69--74. \MR{0389977}

\bibitem{MR0349816}
L.~L. Avramov, \emph{The {H}opf algebra of a local ring}, Izv. Akad. Nauk SSSR
  Ser. Mat. \textbf{38} (1974), 253--277. \MR{0349816}

\bibitem{MR485906}
\bysame, \emph{Small homomorphisms of local rings}, J. Algebra \textbf{50}
  (1978), no.~2, 400--453. \MR{485906}

\bibitem{MR601460}
\bysame, \emph{Obstructions to the existence of multiplicative structures on
  minimal free resolutions}, Amer. J. Math. \textbf{103} (1981), no.~1, 1--31.
  \MR{601460}

\bibitem{MR749041}
\bysame, \emph{Local algebra and rational homotopy}, Algebraic homotopy and
  local algebra ({L}uminy, 1982), Ast\'{e}risque, vol. 113, Soc. Math. France,
  Paris, 1984, pp.~15--43. \MR{749041}

\bibitem{MR846439}
\bysame, \emph{Golod homomorphisms}, Algebra, algebraic topology and their
  interactions ({S}tockholm, 1983), Lecture Notes in Math., vol. 1183,
  Springer, Berlin, 1986, pp.~59--78. \MR{846439}

\bibitem{avramov:vpd}
\bysame, \emph{Modules of finite virtual projective dimension}, Invent. Math.
  \textbf{96} (1989), no.~1, 71--101. \MR{981738}

\bibitem{avramov:ifr}
\bysame, \emph{Infinite free resolutions}, Six lectures on commutative algebra
  (Bellaterra, 1996), Progr. Math., vol. 166, Birkh\"auser, Basel, 1998,
  pp.~1--118. \MR{99m:13022}

\bibitem{avramov:lcih}
\bysame, \emph{Locally complete intersection homomorphisms and a conjecture of
  {Q}uillen on the vanishing of cotangent homology}, Ann. of Math. (2)
  \textbf{150} (1999), no.~2, 455--487. \MR{1726700 (2001a:13024)}

\bibitem{avramov:cslrec3}
\bysame, \emph{A cohomological study of local rings of embedding codepth 3}, J.
  Pure Appl. Algebra \textbf{216} (2012), no.~11, 2489--2506. \MR{2927181}

\bibitem{avramov:svcci}
L.~L. Avramov and R.-O. Buchweitz, \emph{Support varieties and cohomology over
  complete intersections}, Invent. Math. \textbf{142} (2000), no.~2, 285--318.
  \MR{1794064 (2001j:13017)}

\bibitem{avramov:edcrcvct}
L.~L. Avramov, R.-O.\ Buchweitz, and L.~M. {\c{S}}ega, \emph{Extensions of a
  dualizing complex by its ring: commutative versions of a conjecture of
  {T}achikawa}, J. Pure Appl. Algebra \textbf{201} (2005), no.~1-3, 218--239.
  \MR{2158756 (2006e:13012)}

\bibitem{avramov:glh}
L.~L. Avramov and H.-B.\ Foxby, \emph{Gorenstein local homomorphisms}, Bull.
  Amer. Math. Soc. (N.S.) \textbf{23} (1990), no.~1, 145--150. \MR{1020605
  (90k:13009)}

\bibitem{avramov:lgh}
\bysame, \emph{Locally {G}orenstein homomorphisms}, Amer. J. Math. \textbf{114}
  (1992), no.~5, 1007--1047. \MR{1183530 (93i:13019)}

\bibitem{avramov:rhafgd}
\bysame, \emph{Ring homomorphisms and finite {G}orenstein dimension}, Proc.
  London Math. Soc. (3) \textbf{75} (1997), no.~2, 241--270. \MR{98d:13014}

\bibitem{avramov:solh}
L.~L. Avramov, H.-B.\ Foxby, and B.\ Herzog, \emph{Structure of local
  homomorphisms}, J. Algebra \textbf{164} (1994), 124--145. \MR{95f:13029}

\bibitem{avramov:bsolrhoffd}
L.~L. Avramov, H.-B.Foxby, and J.\ Lescot, \emph{Bass series of local ring
  homomorphisms of finite flat dimension}, Trans. Amer. Math. Soc. \textbf{335}
  (1993), no.~2, 497--523. \MR{93d:13026}

\bibitem{avramov:tlg}
L.~L. Avramov and S.\ Halperin, \emph{Through the looking glass: a dictionary
  between rational homotopy theory and local algebra}, Algebra, algebraic
  topology and their interactions ({S}tockholm, 1983), Lecture Notes in Math.,
  vol. 1183, Springer, Berlin, 1986, pp.~1--27. \MR{846435 (87k:55015)}

\bibitem{MR1269426}
L.~L. Avramov and J.~Herzog, \emph{Jacobian criteria for complete
  intersections. {T}he graded case}, Invent. Math. \textbf{117} (1994), no.~1,
  75--88. \MR{1269426}

\bibitem{MR1779800}
L.~L. Avramov and S.~B. Iyengar, \emph{Finite generation of {H}ochschild
  homology algebras}, Invent. Math. \textbf{140} (2000), no.~1, 143--170.
  \MR{1779800}

\bibitem{MR1977825}
\bysame, \emph{Andr\'{e}-{Q}uillen homology of algebra retracts}, Ann. Sci.
  \'{E}cole Norm. Sup. (4) \textbf{36} (2003), no.~3, 431--462. \MR{1977825}

\bibitem{MR2565548}
L.~L. Avramov, S.~B. Iyengar, J.~Lipman, and S.~Nayak, \emph{Reduction of
  derived {H}ochschild functors over commutative algebras and schemes}, Adv.
  Math. \textbf{223} (2010), no.~2, 735--772. \MR{2565548}

\bibitem{avramov:phcnr}
L.~L. Avramov, S.~B. Iyengar, S.~Nasseh, and S.~Sather-Wagstaff,
  \emph{Persistence of homology over commutative noetherian rings}, preprint
  (2020), \texttt{arxiv:2005.10808}.

\bibitem{avramov:htecdga}
\bysame, \emph{Homology over trivial extensions of commutative {DG} algebras},
  Comm. Algebra \textbf{47} (2019), no.~6, 2341--2356. \MR{3957101}

\bibitem{avramov:psmlrsec}
L.~L. Avramov, A.~R. Kustin, and M.~Miller, \emph{Poincar\'e series of modules
  over local rings of small embedding codepth or small linking number}, J.
  Algebra \textbf{118} (1988), no.~1, 162--204. \MR{961334 (89k:13013)}

\bibitem{MR1618363}
D.~Bayer, I.~Peeva, and B.~Sturmfels, \emph{Monomial resolutions}, Math. Res.
  Lett. \textbf{5} (1998), no.~1-2, 31--46. \MR{1618363}

\bibitem{borna:vhgbscm}
K.~Borna, S.\ Sather-Wagstaff, and S.\ Yassemi, \emph{Rings that are
  homologically of minimal multiplicity}, Comm. Algebra \textbf{39} (2011),
  no.~3, 782--807. \MR{2782564}

\bibitem{briggs2020vasconcelos}
B.~Briggs, \emph{Vasconcelos' conjecture on the conormal module}, preprint
  (2020), \texttt{arXiv:2006.04247}.

\bibitem{BenSri}
B.~Briggs and S.~B. Iyengar, \emph{Rigidity properties of the cotangent
  complex}, preprint (2020), \texttt{arXiv:2010.13314}.

\bibitem{BenBriggs}
B.~Briggs, S.~B. Iyengar, J.~C. Letz, and J.~Pollitz, \emph{Locally complete
  intersection maps and the proxy small property}, preprint (2020),
  \texttt{arXiv:2007.08562}.

\bibitem{brion:rq}
M.~Brion, \emph{Representations of quivers}, Geometric methods in
  representation theory. {I}, S\'emin. Congr., vol.~24, Soc. Math. France,
  Paris, 2012, pp.~103--144. \MR{3202702}

\bibitem{bruns:cmr}
W.\ Bruns and J.\ Herzog, \emph{Cohen-{M}acaulay rings}, revised ed., Studies
  in Advanced Mathematics, vol.~39, University Press, Cambridge, 1998.
  \MR{1251956 (95h:13020)}

\bibitem{buchsbaum:asffr}
D.~A. Buchsbaum and D.~Eisenbud, \emph{Algebra structures for finite free
  resolutions, and some structure theorems for ideals of codimension {$3$}},
  Amer. J. Math. \textbf{99} (1977), no.~3, 447--485. \MR{0453723 (56 \#11983)}

\bibitem{MR0077480}
H.~Cartan and S.~Eilenberg, \emph{Homological algebra}, Princeton University
  Press, Princeton, N. J., 1956. \MR{0077480}

\bibitem{CDtest}
O.~Celikbas, H.~Dao, and R.~Takahashi, \emph{Modules that detect finite
  homological dimensions}, Kyoto J. Math. \textbf{54} (2014), no.~2, 295--310.
  \MR{3215569}

\bibitem{celikbas:tgp}
O.~Celikbas and S.~Sather-Wagstaff, \emph{Testing for the {G}orenstein
  property}, Collect. Math. \textbf{67} (2016), no.~3, 555--568. \MR{3536062}

\bibitem{christensen:cmafsdm}
L.~W. Christensen and S.\ Sather-Wagstaff, \emph{A {C}ohen-{M}acaulay algebra
  has only finitely many semidualizing modules}, Math. Proc. Cambridge Philos.
  Soc. \textbf{145} (2008), no.~3, 601--603. \MR{2464778}

\bibitem{christensen:dvke}
\bysame, \emph{Descent via {K}oszul extensions}, J. Algebra \textbf{322}
  (2009), no.~9, 3026--3046. \MR{2567408}

\bibitem{christensen:gmirlr}
L.~W. Christensen, J.\ Striuli, and O.\ Veliche, \emph{Growth in the minimal
  injective resolution of a local ring}, J. Lond. Math. Soc. (2) \textbf{81}
  (2010), no.~1, 24--44. \MR{2580452}

\bibitem{MR3160716}
L.~W. Christensen and O.~Veliche, \emph{Local rings of embedding codepth 3.
  {E}xamples}, Algebr. Represent. Theory \textbf{17} (2014), no.~1, 121--135.
  \MR{3160716}

\bibitem{MR4038053}
L.~W. Christensen, O.~Veliche, and J.~Weyman, \emph{Trimming a {G}orenstein
  ideal}, J. Commut. Algebra \textbf{11} (2019), no.~3, 325--339. \MR{4038053}

\bibitem{MR1803231}
N.~Dupont and M.~Vigu\'{e}-Poirrier, \emph{Finiteness conditions for
  {H}ochschild homology algebra and free loop space cohomology algebra},
  $K$-Theory \textbf{21} (2000), no.~3, 293--300. \MR{1803231}

\bibitem{MR2200850}
W.~G. Dwyer, J.~P.~C. Greenlees, and S.~Iyengar, \emph{Duality in algebra and
  topology}, Adv. Math. \textbf{200} (2006), no.~2, 357--402. \MR{2200850}

\bibitem{MR2225632}
\bysame, \emph{Finiteness in derived categories of local rings}, Comment. Math.
  Helv. \textbf{81} (2006), no.~2, 383--432. \MR{2225632}

\bibitem{eisenbud:ca}
D.~Eisenbud, \emph{Commutative algebra}, Graduate Texts in Mathematics, vol.
  150, Springer-Verlag, New York, 1995, With a view toward algebraic geometry.
  \MR{1322960 (97a:13001)}

\bibitem{MR1037391}
S.~Eliahou and M.~Kervaire, \emph{Minimal resolutions of some monomial ideals},
  J. Algebra \textbf{129} (1990).

\bibitem{MR664027}
Y.~F\'{e}lix and S.~Halperin, \emph{Rational {LS} category and its
  applications}, Trans. Amer. Math. Soc. \textbf{273} (1982), no.~1, 1--38.
  \MR{664027}

\bibitem{MR219546}
D.~Ferrand, \emph{Suite r\'{e}guli\`ere et intersection compl\`ete}, C. R.
  Acad. Sci. Paris S\'{e}r. A-B \textbf{264} (1967), A427--A428. \MR{219546}

\bibitem{foxby:gmarm}
H.-B.\ Foxby, \emph{Gorenstein modules and related modules}, Math. Scand.
  \textbf{31} (1972), 267--284 (1973). \MR{48 \#6094}

\bibitem{frankild:ddgmgdga}
A.~Frankild, S.~Iyengar, and P.~J{\o}rgensen, \emph{Dualizing differential
  graded modules and {G}orenstein differential graded algebras}, J. London
  Math. Soc. (2) \textbf{68} (2003), no.~2, 288--306. \MR{1994683
  (2004f:16013)}

\bibitem{frankild:gdga}
A.~Frankild and P.~J{\o}rgensen, \emph{Gorenstein differential graded
  algebras}, Israel J. Math. \textbf{135} (2003), 327--353. \MR{1997049
  (2005d:16018)}

\bibitem{gabriel:frto}
P.~Gabriel, \emph{Finite representation type is open}, Proceedings of the
  {I}nternational {C}onference on {R}epresentations of {A}lgebras ({C}arleton
  {U}niv., {O}ttawa, {O}nt., 1974), {P}aper {N}o. 10 (Ottawa, Ont.), Carleton
  Univ., 1974, pp.~23 pp. Carleton Math. Lecture Notes, No. 9. \MR{0376769 (51
  \#12944)}

\bibitem{Hugh}
H.~Geller, \emph{{DG} algebra resolutions of fiber products}, in preparation.

\bibitem{gerko:sdc}
A.~A. Gerko, \emph{On the structure of the set of semidualizing complexes},
  Illinois J. Math. \textbf{48} (2004), no.~3, 965--976. \MR{2114263}

\bibitem{golod:hslr}
E.~S. Golod, \emph{Homologies of some local rings}, Dokl. Akad. Nauk SSSR
  \textbf{144} (1962), 479--482. \MR{0138667 (25 \#2110)}

\bibitem{grothendieck:ega3-1}
A.\ Grothendieck, \emph{\'{E}l\'ements de g\'eom\'etrie alg\'ebrique. {III}.
  \'{E}tude cohomologique des faisceaux coh\'erents. {I}}, Inst. Hautes
  \'Etudes Sci. Publ. Math. (1961), no.~11, 167. \MR{0217085 (36 \#177c)}

\bibitem{grothendieck:ega4-2}
\bysame, \emph{\'{E}l\'ements de g\'eom\'etrie alg\'ebrique. {IV}. \'{E}tude
  locale des sch\'emas et des morphismes de sch\'emas. {II}}, Inst. Hautes
  \'Etudes Sci. Publ. Math. (1965), no.~24, 231. \MR{0199181 (33 \#7330)}

\bibitem{MR364232}
T.~H. Gulliksen, \emph{A change of ring theorem with applications to
  {P}oincar\'{e} series and intersection multiplicity}, Math. Scand.
  \textbf{34} (1974), 167--183. \MR{364232}

\bibitem{MR0262227}
T.~H. Gulliksen and G.~Levin, \emph{Homology of local rings}, Queen's Paper in
  Pure and Applied Mathematics, No. 20, Queen's University, Kingston, Ont.,
  1969. \MR{0262227}

\bibitem{happel:sm}
D.~Happel, \emph{Selforthogonal modules}, Abelian groups and modules (Padova,
  1994), Math. Appl., vol. 343, Kluwer Acad. Publ., Dordrecht, 1995,
  pp.~257--276. \MR{1378204 (97d:16016)}

\bibitem{herzog:kadla}
J.~Herzog, \emph{Komplexe, aufl\"{o}sungen, und dualit\"{a}t in der lokalen
  algebra}, Habilitationsschrift, Regensburg, 1973.

\bibitem{herzogv}
\bysame, \emph{Homological properties of the module of differentials}, Atas da
  $6^a$ Escola de \'{A}lgebra (Recife) Cole\c{c}, Atas Soc Brasil. Mat
  \textbf{14} (1981), 35--64.

\bibitem{MR142598}
G.~Hochschild, B.~Kostant, and A.~Rosenberg, \emph{Differential forms on
  regular affine algebras}, Trans. Amer. Math. Soc. \textbf{102} (1962),
  383--408. \MR{142598}

\bibitem{huneke:voeatoscmlr}
C.~Huneke, L.~M. \c{S}ega, and A.~N. Vraciu, \emph{Vanishing of {E}xt and {T}or
  over some {C}ohen-{M}acaulay local rings}, Illinois J. Math. \textbf{48}
  (2004), no.~1, 295--317. \MR{2048226}

\bibitem{huneke:vtci}
C.~Huneke, D.~A. Jorgensen, and R.~Wiegand, \emph{Vanishing theorems for
  complete intersections}, J. Algebra \textbf{238} (2001), no.~2, 684--702.
  \MR{1823780 (2002h:13025)}

\bibitem{MR2355775}
S.~B. Iyengar, \emph{Andr\'{e}-{Q}uillen homology of commutative algebras},
  Interactions between homotopy theory and algebra, Contemp. Math., vol. 436,
  Amer. Math. Soc., Providence, RI, 2007, pp.~203--234. \MR{2355775}

\bibitem{jorgensen:gabf}
D.~A. Jorgensen, \emph{A generalization of the {A}uslander-{B}uchsbaum
  formula}, J. Pure Appl. Algebra \textbf{144} (1999), no.~2, 145--155.
  \MR{1732626 (2000k:13010)}

\bibitem{jorgensen:fpdve}
\bysame, \emph{Finite projective dimension and the vanishing of {${\rm
  Ext}_R(M,M)$}}, Comm. Algebra \textbf{36} (2008), no.~12, 4461--4471.
  \MR{2473341 (2009i:13026)}

\bibitem{jorgensen:gbscm}
D.~A. Jorgensen and G.~J. Leuschke, \emph{On the growth of the {B}etti sequence
  of the canonical module}, Math. Z. \textbf{256} (2007), no.~3, 647--659.
  \MR{2299575 (2008a:13018)}

\bibitem{jorgensen:nccgr}
D.~A. Jorgensen and L.~M. {\c{S}}ega, \emph{Nonvanishing cohomology and classes
  of {G}orenstein rings}, Adv. Math. \textbf{188} (2004), no.~2, 470--490.
  \MR{2087235 (2005f:13017)}

\bibitem{jorgensen:itrcm}
\bysame, \emph{Independence of the total reflexivity conditions for modules},
  Algebr. Represent. Theory \textbf{9} (2006), no.~2, 217--226. \MR{2238367
  (2007c:13022)}

\bibitem{MR3862675}
L.~Katth\"{a}n, \emph{The structure of {DGA} resolutions of monomial ideals},
  J. Pure Appl. Algebra \textbf{223} (2019), no.~3, 1227--1245. \MR{3862675}

\bibitem{kustin:gacfct}
A.~R. Kustin, \emph{Gorenstein algebras of codimension four and characteristic
  two}, Comm. Algebra \textbf{15} (1987), no.~11, 2417--2429. \MR{912779
  (88j:13020)}

\bibitem{MR1132435}
\bysame, \emph{Classification of the {T}or-algebras of codimension four almost
  complete intersections}, Trans. Amer. Math. Soc. \textbf{339} (1993), no.~1,
  61--85. \MR{1132435}

\bibitem{MR1292771}
\bysame, \emph{The minimal resolution of a codimension four almost complete
  intersection is a {DG}-algebra}, J. Algebra \textbf{168} (1994), no.~2,
  371--399. \MR{1292771}

\bibitem{kustin:asmrgrecf}
A.~R. Kustin and M.~Miller, \emph{Algebra structures on minimal resolutions of
  {G}orenstein rings of embedding codimension four}, Math. Z. \textbf{173}
  (1980), no.~2, 171--184. \MR{583384 (81j:13013)}

\bibitem{MR713381}
\bysame, \emph{Multiplicative structure on resolutions of algebras defined by
  {H}erzog ideals}, J. London Math. Soc. (2) \textbf{28} (1983), no.~2,
  247--260. \MR{713381}

\bibitem{MR1295961}
A.~R. Kustin and S.~M. Palmer~Slattery, \emph{The {P}oincar\'{e} series of
  every finitely generated module over a codimension four almost complete
  intersection is a rational function}, J. Pure Appl. Algebra \textbf{95}
  (1994), no.~3, 271--295. \MR{1295961}

\bibitem{Levin}
G.~Levin, \emph{Finitely generated {${\rm Ext}$} algebras}, Math. Scand.
  \textbf{49} (1981), no.~2, 161--180 (1982). \MR{661889}

\bibitem{MR1217970}
J.-L. Loday, \emph{Cyclic homology}, Grundlehren der Mathematischen
  Wissenschaften [Fundamental Principles of Mathematical Sciences], vol. 301,
  Springer-Verlag, Berlin, 1992, Appendix E by Mar\'{\i}a O. Ronco.
  \MR{1217970}

\bibitem{Todd}
T.~Morra, \emph{{DG} algebra structures on resolutions of stanley-reisner
  ideals of certain simplicial spheres}, in preparation.

\bibitem{MR104718}
T.~Nakayama, \emph{On algebras with complete homology}, Abh. Math. Sem. Univ.
  Hamburg \textbf{22} (1958), 300--307. \MR{104718}

\bibitem{nasseh:ldgm}
S.~Nasseh and S.~Sather-Wagstaff, \emph{Liftings and quasi-liftings of {DG}
  modules}, J. Algebra \textbf{373} (2013), 162--182. \MR{2995021}

\bibitem{nasseh:egdgm}
\bysame, \emph{Extension groups for {DG} modules}, Comm. Algebra \textbf{45}
  (2017), no.~10, 4466--4476. \MR{3640821}

\bibitem{nasseh:gart}
\bysame, \emph{Geometric aspects of representation theory for {DG} algebras:
  answering a question of {V}asconcelos}, J. Lond. Math. Soc. (2) \textbf{96}
  (2017), no.~1, 271--292. \MR{3687949}

\bibitem{nasseh:vetfp}
\bysame, \emph{Vanishing of {E}xt and {T}or over fiber products}, Proc. Amer.
  Math. Soc. \textbf{145} (2017), no.~11, 4661--4674. \MR{3691985}

\bibitem{nasseh:ahplrdmi}
S.~Nasseh, S.~Sather-Wagstaff, R.~Takahashi, and K.~VandeBogert,
  \emph{Applications and homological properties of local rings with
  decomposable maximal ideals}, J. Pure Appl. Algebra \textbf{223} (2019),
  no.~3, 1272--1287. \MR{3862678}

\bibitem{nasseh:lrqdmi}
S.~Nasseh and R.~Takahashi, \emph{Local rings with quasi-decomposable maximal
  ideal}, Math. Proc. Cambridge Philos. Soc. \textbf{168} (2020), no.~2,
  305--322. \MR{4064107}

\bibitem{nasseh:oeire}
S.~Nasseh and Y.~Yoshino, \emph{On {E}xt-indices of ring extensions}, J. Pure
  Appl. Algebra \textbf{213} (2009), no.~7, 1216--1223. \MR{2497570
  (2010f:13016)}

\bibitem{MR3774891}
\bysame, \emph{Weak liftings of {DG} modules}, J. Algebra \textbf{502} (2018),
  233--248. \MR{3774891}

\bibitem{MR4152766}
M.~Ono and Y.~Yoshino, \emph{A lifting problem for {DG} modules}, J. Algebra
  \textbf{566} (2021), 342--360. \MR{4152766}

\bibitem{MR1407879}
I.~Peeva, \emph{{$0$}-{B}orel fixed ideals}, J. Algebra \textbf{184} (1996),
  no.~3, 945--984. \MR{1407879}

\bibitem{Planas}
F.~Planas-Vilanova, \emph{On the vanishing and non-rigidity of the
  {A}ndr\'{e}-{Q}uillen (co)homology}, J. Pure Appl. Algebra \textbf{120}
  (1997), no.~1, 67--75. \MR{1466098}

\bibitem{Josh2}
J.~Pollitz, \emph{Cohomological supports over derived complete intersections
  and local rings}, preprint (2019), \texttt{arXiv:1912.12009}.

\bibitem{MR3988642}
J.~Pollitz, \emph{The derived category of a locally complete intersection
  ring}, Adv. Math. \textbf{354} (2019), 106752, 18. \MR{3988642}

\bibitem{MR0257068}
D.~Quillen, \emph{On the (co-) homology of commutative rings}, Applications of
  {C}ategorical {A}lgebra ({P}roc. {S}ympos. {P}ure {M}ath., {V}ol. {XVII},
  {N}ew {Y}ork, 1968), Amer. Math. Soc., Providence, R.I., 1970, pp.~65--87.
  \MR{0257068}

\bibitem{test}
S.~Sather-Wagstaff, \emph{Ascent properties for test modules}, preprint (2019),
  \texttt{arXiv:1911.07708}.

\bibitem{sather:cidfc}
\bysame, \emph{Complete intersection dimensions and {F}oxby classes}, J. Pure
  Appl. Algebra \textbf{212} (2008), no.~12, 2594--2611. \MR{2452313
  (2009h:13015)}

\bibitem{sather:bnsc}
\bysame, \emph{Bass numbers and semidualizing complexes}, Commutative algebra
  and its applications, Walter de Gruyter, Berlin, 2009, pp.~349--381.
  \MR{2640315}

\bibitem{MR1974627}
L.~M. {\c{S}}ega, \emph{Vanishing of cohomology over {G}orenstein rings of
  small codimension}, Proc. Amer. Math. Soc. \textbf{131} (2003), no.~8,
  2313--2323. \MR{1974627}

\bibitem{sega:stfcar}
\bysame, \emph{Self-tests for freeness over commutative {A}rtinian rings}, J.
  Pure Appl. Algebra \textbf{215} (2011), no.~6, 1263--1269. \MR{2769231}

\bibitem{serre:sldhdaedmn}
J.-P. Serre, \emph{Sur la dimension homologique des anneaux et des modules
  noeth\'eriens}, Proceedings of the international symposium on algebraic
  number theory, Tokyo \& Nikko, 1955 (Tokyo), Science Council of Japan, 1956,
  pp.~175--189. \MR{19,119a}

\bibitem{skoldberg1}
E.~Sk\"{o}ldberg, \emph{Resolutions of modules with initially linear syzygies},
  preprint (2011), \texttt{arXiv:1106.1913v2}.

\bibitem{MR0349740}
H.~Tachikawa, \emph{Quasi-{F}robenius rings and generalizations. {${\rm QF}-3$}
  and {${\rm QF}-1$} rings}, Lecture Notes in Mathematics, Vol. 351,
  Springer-Verlag, Berlin-New York, 1973, Notes by Claus Michael Ringel.
  \MR{0349740}

\bibitem{tate:hnrlr}
J.~Tate, \emph{Homology of {N}oetherian rings and local rings}, Illinois J.
  Math. \textbf{1} (1957), 14--27. \MR{0086072 (19,119b)}

\bibitem{Ehsan}
E.~Tavanfar, \emph{A trilogy, given by complete tensor product of complete
  rings over the coefficient ring}, preprint (2020), \texttt{arXiv:1911.11290}.

\bibitem{taylor:igmrs}
D.~K. Taylor, \emph{Ideals generated by monomials in an $r$-sequence}, ProQuest
  LLC, Ann Arbor, MI, 1966, Thesis (Ph.D.)--The University of Chicago.
  \MR{2611561}

\bibitem{kellerV}
K.~VandeBogert, \emph{Resolution and {T}or algebra structures of grade $3$
  ideals defining compressed rings}, preprint (2020),
  \texttt{arXiv:2004.06691}.

\bibitem{MR213345}
W.~V. Vasconcelos, \emph{Ideals generated by {$R$}-sequences}, J. Algebra
  \textbf{6} (1967), 309--316. \MR{213345}

\bibitem{vasconcelos:dtmc}
\bysame, \emph{Divisor theory in module categories}, North-Holland Publishing
  Co., Amsterdam, 1974, North-Holland Mathematics Studies, No. 14, Notas de
  Matem\'atica No. 53. [Notes on Mathematics, No. 53]. \MR{0498530 (58
  \#16637)}

\bibitem{MR508082}
\bysame, \emph{On the homology of {$I/I^{2}$}}, Comm. Algebra \textbf{6}
  (1978), no.~17, 1801--1809. \MR{508082}

\bibitem{MR814190}
\bysame, \emph{The complete intersection locus of certain ideals}, J. Pure
  Appl. Algebra \textbf{38} (1985), no.~2-3, 367--378. \MR{814190}

\bibitem{MR1245091}
M.~Vigu\'{e}-Poirrier, \emph{Crit\`eres de nullit\'{e} pour l'homologie des
  alg\`ebres gradu\'{e}es}, C. R. Acad. Sci. Paris S\'{e}r. I Math.
  \textbf{317} (1993), no.~7, 647--649. \MR{1245091}

\bibitem{voigt:idteag}
D.~Voigt, \emph{Induzierte {D}arstellungen in der {T}heorie der endlichen,
  algebraischen {G}ruppen}, Lecture Notes in Mathematics, Vol. 592,
  Springer-Verlag, Berlin, 1977, Mit einer englischen Einf{\"u}hrung.
  \MR{0486168 (58 \#5949)}

\bibitem{MR319985}
J.~Watanabe, \emph{A note on {G}orenstein rings of embedding codimension
  three}, Nagoya Math. J. \textbf{50} (1973), 227--232. \MR{319985}

\bibitem{weyman:sfrl3}
J.~Weyman, \emph{On the structure of free resolutions of length {$3$}}, J.
  Algebra \textbf{126} (1989), no.~1, 1--33. \MR{1023284 (91g:13019)}

\bibitem{yoshino}
Y.\ Yoshino, \emph{The theory of {L}-complexes and weak liftings of complexes},
  J. Algebra \textbf{188} (1997), no.~1, 144--183. \MR{98i:13024}

\end{thebibliography}
\providecommand{\bysame}{\leavevmode\hbox to3em{\hrulefill}\thinspace}
\providecommand{\MR}{\relax\ifhmode\unskip\space\fi MR }
\providecommand{\MRhref}[2]{%
  \href{http://www.ams.org/mathscinet-getitem?mr=#1}{#2}
}
\providecommand{\href}[2]{#2}

\end{document}